\documentclass[MSNbibl,nameyear,seceqn,dvips]{arxstspdf}
\usepackage{mathbh}
\usepackage{graphicx}
\usepackage{flushend}
\usepackage{stfloats}
\usepackage{mathrsfs}

\volume{28}
\issue{4}
\pubyear{2013}
\firstpage{487}
\lastpage{509}
\doi{10.1214/13-STS442} 

\makeatletter

\newcommand{\citecs}[1]{\citeauthor{#1}, \citeyear{#1}}

\newcommand{\iint}{\int\!\!\int}

\newcommand{\rright}{\right}
\newcommand{\lleft}{\left}

\newtheorem{theorem}{Theorem}[section]
\newtheorem{lemma}{Lemma}[section]
\newtheorem{cor}{Corollary}[section]

\newproclaim{remark}{Remark}[section]
\newproclaim{definition}{Definition}[section]

\newcommand{\cal}{\mathcal}
\renewcommand{\citep}[1]{(\citeauthor{#1}, \citeyear{#1})}

\newcommand{\sign}{\operatorname{sign}}
\newcommand{\sud}{\,\mathrm{d}}
\newcommand{\ud}{\mathrm{d}}
\newcommand{\deriv}[3][]{\frac{\ud^{#1} #2}{\ud{#3}^{#1}}}
\newcommand{\pderiv}[3][]{\frac{\partial^{#1} #2}{\partial{#3}^{#1}}}
\renewcommand{\div}[1][]{\nabla_{#1}\cdot}

\renewcommand{\P}{\mathbb P}
\newcommand{\EXP}{\mathbb{E}}
\newcommand{\IND}[1]{{\mathbh{1}_{#1}}}

\newcommand{\vp}{\varphi}
\newcommand{\eps}{\varepsilon}

\newcommand{\Z}{\mathbb{Z}}
\newcommand{\R}{\mathbb{R}}
\newcommand{\Ab}{\bolds{\mathcal{A}}}
\newcommand{\A}{\mathcal{A}}
\newcommand{\E}{\mathcal{E}}
\newcommand{\Dom}{\operatorname{Dom}}

\newcommand{\Lp}[1]{L^{#1}}
\newcommand{\ab}{\bolds{\alpha}}

\newcommand{\x}{{\mathbf x}}

\newcommand{\Db}{\bar{D}}
\newcommand{\Cb}{\bar{C}}

\makeatother

\begin{document}
\begin{frontmatter}

\title{Advection--Dispersion Across Interfaces}
\runtitle{Advection--Dispersion}

\begin{aug}
\author[a]{\fnms{Jorge M.} \snm{Ramirez}\ead[label=e1]{jmramirezo@unalmed.edu.co}},
\author[b]{\fnms{Enrique A.} \snm{Thomann}\ead[label=e2]{thomann@math.oregonstate.edu}} 
\and
\author[b]{\fnms{Edward C.} \snm{Waymire}\corref{}\ead[label=e3]{waymire@math.oregonstate.edu}}
\runauthor{J. M. Ramirez, E. A. Thomann and E. C. Waymire}

\affiliation{Universidad Nacional de Colombia,
Sede Medell\'{i}n and Oregon State University}

\address[a]{Jorge~M. Ramirez is Associate Professor of Mathematics,
Escuela de Matem\'{a}ticas, Universidad Nacional de Colombia
Sede Medell\'{i}n, Medell\'{i}n, Colombia
\printead{e1}.}
\address[b]{Enrique~A. Thomann and Edward~C. Waymire
are Professors of Mathematics, Department of Mathematics,
Oregon State University, Corvallis, Oregon, USA
\printead{e2};\break \printead*{e3}.}
\end{aug}

%
\begin{abstract}
This article concerns a systemic manifestation of small scale
interfacial heterogeneities in large scale quantities of interest to a
variety of diverse applications spanning the earth, biological and
ecological sciences. Beginning with formulations in terms of partial
differential equations governing the conservative, advective-dispersive
transport of mass concentrations in divergence form, the specific
interfacial heterogeneities are introduced in terms of (spatial)
discontinuities in the diffusion coefficient across a lower-dimensional
hypersurface. A pathway to an equivalent stochastic formulation is
then developed with special attention to the interfacial effects in
various functionals such as first passage times, occupation times and
local times. That an appreciable theory is achievable within a
framework of applications involving one-dimensional models having
piecewise constant coefficients greatly facilitates our goal of a
gentle introduction to some rather dramatic mathematical consequences
of interfacial effects that can be used to predict structure and to
inform modeling.
\end{abstract}

%

\end{frontmatter}

\section{Introduction}\label{sec1}

To set the perspective of this article, let us first consider classic
advection--dispersion
phenomena in $ \mathbb{R}^k$ of a concentration of particles immersed
in a fluid
as described by the following\vspace*{1pt} partial differential equation for $
{\mathbf x}\in\mathbb{R}^k$ and $ t \ge0$:
%
%
\begin{eqnarray}
\label{adpde}
\pderiv{u} {t} &=& \frac{1}{2}\nabla\cdot\bigl( {\mathbf D}({\mathbf
x})\nabla u(t,{\mathbf x})\bigr) \nonumber\\
&&{}- \nabla\cdot\bigl({\mathbf v}({\mathbf x})u(t,{\mathbf x})
\bigr),\\
&&\eqntext{u\bigl(0^+,{\mathbf x}\bigr) = u_0({\mathbf x}).}%
\end{eqnarray}
In particular, assume that the coefficients ${\mathbf D}$ and ${\mathbf v}$ are
smooth (matrix/vector-valued) functions\footnote{Throughout this
article we restrict attention to time homogeneous equations.} on
$\mathbb{R}^k$,
$\nabla= \sum_{j=1}^k\pderiv{}{x_j}$. Such\vspace*{1pt} an equation describes the evolution
of an initial (scalar) mass concentration
$u_0$ evolving at a temporal rate assumed to be locally controlled
by spatial fluxes $\frac{1}{2}{\mathbf D}({\mathbf x})\nabla u(t,{\mathbf x}) -
{\mathbf v}({\mathbf x})u(t,{\mathbf x})$. The first term expresses
Fick's law of flux as being proportional to the concentration
gradient, and the second term being
the advection of mass by fluid velocity. Many physical as well as
biological/ecological problems take
this form, perhaps on a spatial domain ${\mathbf G}\subseteq\mathbb
{R}^k$, with
appropriate boundary conditions, such as Dirichlet or Neumann boundary
conditions. The success of nineteenth- and twentieth-century
mathematical developments in\break  analysis, geometry and numerical
computation continues to guide research into the 21st century. In
particular, this
single linear partial differential equation has inspired a body of
mathematical research that
is likely unmatched in diversity and scope.

The pervasive role of ({\ref{adpde}) in the development of probability
and statistical theory is therefore
no surprise. The recognition of the fundamental role of standard
Brownian motion ${\mathbf B}
= \{{\mathbf B}(t)\dvtx  t\ge0\}$ and the corresponding
It\^o's stochastic calculus opened the door to a more natural
reformulation of advective-dispersive
phenomena in terms of the stochastic differential equation
%
%
\begin{eqnarray}
\label{adsde} \ud{\mathbf X}(t) = \tilde{\mathbf v}\bigl({\mathbf X}(t)\bigr) \sud t +
\sqrt{{\mathbf D}\bigl({\mathbf X}(t)\bigr)} \sud{\mathbf B}(t),\nonumber\\[-8pt]\\[-8pt]
\eqntext{t > 0, {\mathbf X}(0) = {\mathbf x},}
\end{eqnarray}
relating the conditional distribution $p(t,{\mathbf x},\ud{\mathbf y})$ of
${\mathbf X}(t)$ given ${\mathbf X}(0) = {\mathbf x}$, that is, the transition
probabilities,
to the fundamental solution of (\ref{adpde}) via
the basic semigroup formula
%
%
\begin{equation}
u(t,{\mathbf x}) = \int_{\mathbb R^k}u_0({\mathbf y})p(t,{\mathbf
x},\ud{\mathbf y}).
\end{equation}
In (\ref{adsde}), $\sqrt{{\mathbf D}}$ is the matrix square root of the
molecular dispersion tensor ${\mathbf D}$, which augments the macroscopic advection
${\mathbf v}$ via
%
%
\begin{equation}
\label{microdrift}\quad \tilde{\mathbf v}({\mathbf x}) = \biggl(- \sum
_j \,\frac{1}{2}\pderiv {D_{ji}}
{x_j}(\x) + v_i({\mathbf x}) \biggr)_{1\le i\le k}.
\end{equation}
The Markov process ${\mathbf X}$ so determined becomes the probabilistic
representation of the object of interest in relation to the p.d.e. (\ref
{adpde}), be it physical, biological or perhaps financial.

Not only does a stochastic framework enable new approaches to the
analysis of problems
related to (\ref{adpde}), but it inspires still more diverse ways in
which to model, analyze
and measure naturally occurring phenomena. After all, the coefficients
${\mathbf v}$ and ${\mathbf D}$ now admit a statistical interpretation! The
significance of this fact
was made manifestly clear through the observations and
measurements of \citet{perrin} in his
historic determination of Avogadro's constant, following up Einstein's
1905 theory of the molecular structure of matter.
In addition,
new models that may be a priori less obvious to formulate at the scale
of (\ref{adpde}) emerge to describe phenomena at the scale of particle
trajectories as
observed in certain financial data or biological experiments
(see, e.g., \citecs{decamps2006}; Fagan, Cantrell and Cosner, \citeyear{fagan1}).
Moreover, in the context of particle trajectories, a wide variety of
sample path functionals, such as first passage times, escape and
occupation times, and local times also emerge naturally in both theory
and applications.

From a probabilistic perspective, the smoothness of the coefficients in
(\ref{adpde}) goes a
long way toward the alternative view expressed through (\ref{adsde}) of
particle trajectories being (approximately) shifts and rescalings of a
standard Brownian
motion when observed locally (infinitesimally) in time. In particular,
if the coefficients
are in fact constant, then the solution to (\ref{adsde}) is a Brownian motion
%
%
\begin{equation}
\label{genbm} {\mathbf X}(t) = {\mathbf x} + {\mathbf v}t
+ \sqrt{{\mathbf D}} {\mathbf B}(t),\quad
t\ge0,
\end{equation}
with \textit{drift coefficient} $\tilde{\mathbf v} \equiv{\mathbf v}$,
and \textit{diffusion coeffi-\break cient}~${\mathbf D}$,
whose transition probabilities $p(t,{\mathbf x},{\mathbf y})$,
assuming nonsingularity ($\det{\mathbf D} \neq0$),
provide the fundamental solution to
%
%
\begin{equation}
\label{adpdeconstantcoeff} \pderiv{u} {t} = \frac{1}{2}{\mathbf D}\Delta u + {
\mathbf v}\cdot\nabla u
\end{equation}
with Laplacian $\Delta\equiv\nabla\cdot\nabla= \sum_{j=1}^k \pderiv
[2]{}{x_j}$, that is,
differentiation with respect to the backward variable ${\mathbf x}$. More
generally, assuming sufficient smoothness of coefficients,
(\ref{adpde}) may be directly recast after relabeling ${\mathbf x}$ as
${\mathbf y}$,
in the form of \textit{Kolmogorov's forward equation}, or the \textit{Fokker--Planck equation} as
it is called in the physical sciences,
%
%
\begin{eqnarray}
\label{forwardadpde} \pderiv{u} {t}(t, {\mathbf y}) &=& \frac{1}{2} \sum
_{i,j}\frac{\partial^2
(D_{ij}({\mathbf y})u({t,\mathbf y}))}{\partial y_i\,\partial y_j}\nonumber\\[-8pt]\\[-8pt]
&&{}- \sum
_i \frac
{\partial(\tilde v_i({\mathbf y})u(t,{\mathbf y}))}{\partial
y_i},\nonumber
\end{eqnarray}
where $\tilde{\mathbf v}$ already appeared in
(\ref{microdrift}).
Note that this is merely an equivalent way in which to express
the equation (\ref{adpde}), with the relabeling of variables suggested
by their respective roles as backward
and forward variables
in the transition probabilities $p(t,\break {\mathbf x},d{\mathbf y})$. On the other
hand, \textit{Kolmogorov's backward equation}, with (\ref
{adpdeconstantcoeff}) as a special case, is obtained from
(\ref{forwardadpde}) by integration by parts as the adjoint
%
%
\begin{eqnarray}
\label{backwardadpde} \pderiv{u} {t}(t, {\mathbf x}) &=& \frac{1}{2}\sum
_{i,j}D_{ij}({\mathbf x}) \,\frac{\partial^2 u}{\partial x_i\,\partial x_j}(t,{
\mathbf x}) \nonumber\\[-8pt]\\[-8pt]
&&{}+ \sum_i\tilde v_i({\mathbf x})\,
\frac{\partial u}{\partial x_i}(t,{\mathbf x}).\nonumber
\end{eqnarray}

As will be discussed as the primary point of the present article,
there are phenomena for which the
smoothness of the coefficients is untenable.
The particular ``nonsmoothness''
of focus here can most generally be framed as a discontinuity,
of otherwise\break  (piecewise) continuous coefficients, on a hypersurface of
co-dimension one. This includes discontinuities at \mbox{(0-dimensional)}
points in one dimension or across a curve in two dimensions.

Advection--dispersion was framed in terms of the $k$-dimensional model
(\ref{adpde}) in an effort to frame the \textit{big} problems for
continued research. However, there is much yet to be learned about
interfacial problems in dimensions greater than one.
Perhaps surprisingly, but indeed fortunately, the
applications involving one-dimensional processes are already extensive enough
to provide a rich source of examples with features of both mathematical and
empirical interest, especially as manifested in the behavior of the
functionals noted earlier.
One may also expect that some results in one dimension will at least
partially inform higher-dimensional problems.\looseness=1

Just as standard Brownian motion plays a basic, albeit secondary, role
in constructing the Markov process
$X$ associated with (\ref{adpde}) in the case of smooth coefficients
via (\ref{adsde}), a
class of \textit{skew Brownian motions} will emerge in the construction of
the Markov processes (termed \textit{skew diffusions}) associated with
one-dimensional advection--dispersion across an interface.
\textit{Skew Brownian motion} $B^{(\alpha)} = \{B^{(\alpha)}(t)\dvtx\break   t\ge0\},
0 < \alpha< 1$ is a continuous semimartingale introduced by
\citet{IM1963}. Fundamental papers on skew Brownian motion by
\citet{harrisonshepp}, \citet{walsh1}, \citet{Ouknine}
and \citet{LeGall} are summarized in a mathematically
comprehensive survey article by \citet{lejaysurvey}. Interesting
fresh ideas on some of the foundational questions about skew Brownian
motion continue to emerge, for example, \citet{hairer1},
\citet{Prokaj} and \citet{Fernholz2012fk}.
These provide a number of equivalent ways in which to view skew
Brownian motion on which the present survey article will build.
It is not our intention to provide a mathematically comprehensive
survey of skew Brownian motion.\footnote{The survey article
\citet{lejaysurvey} in fact fills this need quite thoroughly and
is recommended as follow-up to the present article.} Indeed, the
primary focus here is on the Markov process\vadjust{\goodbreak} (skew diffusion)
associated with advection--dispersion across an interface in one dimension.
Our goal is to provide a simple, focused mathematical framework of skew
diffusion in which to then illustrate rather dramatic consequences of
\textit{interfacial effects} pertaining to specific physical and
biological phenomena.

Just as with the case of smooth coefficients,
the analysis, in terms of both partial differential and stochastic
differential equations, the numeric computational schemes, and the
statistical aspects of advection--dispersion across an interface, relate
to each other in an interesting mathematical interplay. Accordingly, as
will be illustrated by examples, empirical observations in this context
often point to new and interesting phenomena amenable to mathematical
explanation or prediction.

Therein lies the overarching goal of this article. Namely, within the
context of the Mathematics of Planet Earth and International Year of
Statistics 2013 initiatives,
we seek to illustrate some of the mathematical structure reflected in
observed and predicted large scale properties of advection--disper\-sion
as a consequence
of locally defined \textit{interfacial discontinuities} of the type
described above. For example, results are described that quantify the
dramatic effect of a (small scale) point discontinuity
on the behavior of occupation times of large scale regions in one dimension.
This is achieved by identification and analysis of the basic Markov
process associated with the given coefficients.
That this is in fact achievable within a framework of one-dimensional
models with piecewise constant coefficients facilitates our goal of a
gentle introduction to results that are also relevant to a diverse
range of applications to be described herein.


The organization of the paper is as follows: in Section~\ref{sec1} skew Brownian
motion and its properties are introduced in a broader context of
dispersion of a solute in the presence of a so-called \textit{conservative
interface condition}, that is, physical skew diffusion. This is
followed by subsequent sections to provide illustrations of some more
general (nonconservative) interface conditions that arise naturally in
the physical and biological sciences, including
free surface heights in ocean upwellings, animal movement models in
ecology and dispersion in a river network. Building on these examples,
a summary of complementary results and open problems inspired by these
examples is provided in the closing section.\looseness=1

\section{One-Dimensional Physical Skew Diffusion and Skew Brownian Motion}\label{sec2}

Building on the theme laid out in the \hyperref[sec1]{Introduction},
the Markov process to be referred to as one-dimen\-sional
\textit{physical skew diffusion} with parameters $D^+>0$, $D^->0$ ($v =
0$) will be defined in relation to the following continuity equation
for a solute immersed in fluid medium separated by
a point interface at the origin:
%
%
\begin{eqnarray}
\label{1dskewpde}
\pderiv{u} {t}(t,x) &=& \frac{1}{2} D(x)\,\pderiv[2]{u}
{x}(t,x)\quad (x\neq 0),
\\
\qquad D^+\,\pderiv{u} {x}\bigl(t,0^+\bigr) &=& D^-\,\pderiv{u} {x}\bigl(t,0^-
\bigr),\nonumber\\
\eqntext{u\bigl(t,0^+\bigr) = u\bigl(t,0^-\bigr), t > 0,}
\end{eqnarray}
where
%
%
\begin{equation}
\label{definitionD} D(x) = \cases{D^+, & if $x > 0$,
\cr
D^-, & if $x\leq0$.}
\end{equation}
The particular interface condition
%
%
\begin{equation}
\label{conserveinterface} D^+\,\pderiv{u} {x}\bigl(t,0^+\bigr) = D^-\,\pderiv{u} {x}
\bigl(t,0^-\bigr)
\end{equation}
ensures that the diffusive \textit{flux} $D(x) \,\pderiv{u}{x}(t,x)$ is
continuous at all $x \in\R$ and for all $t>0$. Moreover, it yields
``conservation of mass'' $\int_{-\infty}^\infty u(t,x)\sud x =\break
\int_{-\infty}^\infty u(0,x)\sud x$, $t>0$, since after integration by
parts, one has
%
%
\begin{equation}
\deriv{} {t}\int_{-\infty}^\infty u(t,x) \sud x = 0.
\end{equation}

In particular, this interface condition makes the spatial operator in
(\ref{1dskewpde})
formally self-adjoint.

The simplest approach to identify the corresponding physical skew
diffusion process
is perhaps by explicitly solving (\ref{1dskewpde}). Indeed, for any
initial condition $u(0,x) = u_0(x)$, equation (\ref{1dskewpde}) has
solution $u(t,y) = \int_{-\infty}^{\infty} p^*(t,x,y) u_0(x) \sud x$
where the fundamental solution $p^*(t,x,y)$ can be simply checked to be
(see \citecs{ramirez2006})
%
%
\begin{eqnarray}
\label{eqskewdifftrans}\qquad
&&p^*(t,x,y) \nonumber\\[-8pt]\\[-8pt]
&&\quad= \cases{\displaystyle\frac{1}{\sqrt{2 \pi D^+ t}}
\biggl[\exp \biggl
\{\frac{-(y - x)^2}{2 D^+
t} \biggr\} \vspace*{2pt}\cr
\quad{} + \displaystyle
\frac{\sqrt{D^+} - \sqrt{D^-}}{\sqrt{D^-}+\sqrt{D^+}} \exp
\biggl\{\frac{-(y + x)^2}{2 D^+t}
\biggr\} \biggr], \cr
\quad\mbox{$x>0, y>0$},
\vspace*{2pt}\cr
\displaystyle\frac{1}{\sqrt{2 \pi D^- t}} \biggl[\exp \biggl\{
\frac{-(y - x)^2}{2 D^-
t} \biggr\} \vspace*{2pt}\cr
\quad{}-\displaystyle\frac{\sqrt{D^+} - \sqrt{D^-}}{\sqrt{D^+}+\sqrt{D^-}} \exp \biggl\{\frac{-(y + x)^2}{2 D^- t}
\biggr\} \biggr], \cr
\quad\mbox{$x<0, y<0$},
\vspace*{2pt}\cr
\displaystyle\frac{2}{\sqrt{D^+}+\sqrt{D^-}} \frac{1}{\sqrt{2\pi t}}\vspace*{2pt}\cr
\quad{}\cdot\displaystyle\exp
\biggl\{ -\frac{(y\sqrt{D^-} - x\sqrt{D^+})^2}{2 D^- D^+ t} \biggr\}, \cr
\quad\mbox{$x \leq0, y\geq0$},
\vspace*{2pt}\cr
\displaystyle\frac{2}{\sqrt{D^+}+\sqrt{D^-}}
\frac{1}{\sqrt{2\pi t}} \vspace*{2pt}\cr
\quad{}\cdot\displaystyle\exp \biggl\{ -\frac{(y\sqrt{D^+} - x\sqrt{D^-})^2}
{2 D^- D^+ t} \biggr\}, \cr
\quad\mbox{$x \geq0, y
\leq0$}.}\nonumber
\end{eqnarray}

See Figure~\ref{Figureptxy}. Observe that while
%
%
\begin{equation}
\label{physym}\qquad p^*(t,x,y) = p^*(t,y,x),\quad x,y\in\R, t > 0,
\end{equation}
there is nonetheless a
``skewness asymmetry'' around the interface exhibited in the calculation
%
%
\begin{equation}
\label{skeweffect} \int_{[0,\infty)}p^*(t,0,y )\sud y =
\frac{\sqrt{D^+}}{\sqrt{D^+}+\sqrt{D^-}}.
\end{equation}

\begin{figure}

\includegraphics{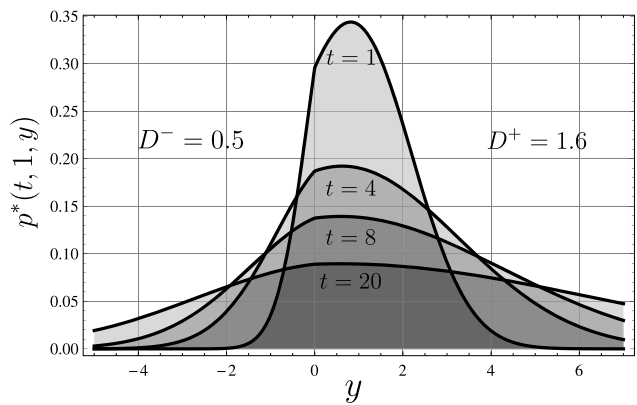}

\caption{Plots of $p^*(t,x,y)$ for $x=1$, several values of $t$k and
fixed $D^+ > D^-$.}\label{Figureptxy}
\end{figure}

To prepare for the definition of skew Brownian motion,
let $B = \{B(t)\dvtx  t\ge0\}$ denote standard Brownian motion started at
$B(0) = 0$,
and let $A = \{A_n\dvtx\break   n=1,2,\ldots\}$ be
an i.i.d. sequence of $\pm1$-valued Ber\-noulli random variables with
$\alpha= P(A_n =1)$, independent of
$B$, defined on a common probability space $(\Omega, {\cal F},P)$.
Since the paths $t\to B(t)$ are continuous,
%
\begin{figure}

\includegraphics{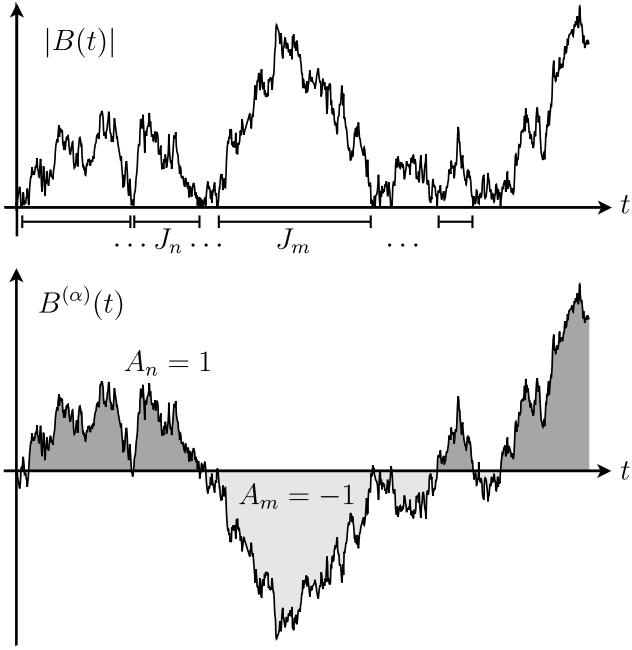}

\caption{Skew Brownian motion construction.} \label{skewbm}\label{FigureSBM}
\end{figure}
the complement to the closed subset $B^{-1}(\{0\})$ is a countable
disjoint union
of (random) open intervals of $[0,\infty)$ enumerated as $J_1,J_2,\ldots\,$;
see Figure~\ref{FigureSBM}.

\begin{definition}[($\alpha$-skew Brownian motion)]
Let $\alpha\in[0,1]$. The stochastic process given by
\[
B^{(\alpha)}(t) = \sum_{n=1}^\infty
A_n \IND{J_n}(t) \bigl|B(t)\bigr|,\quad t\ge0,
\]
is referred to as \textit{skew Brownian motion with transmission parameter
$\alpha$} starting
at $0$.
\end{definition}

\begin{remark}
The cases $\alpha=0,1$ correspond to reflecting Brownian motion
$B^{(0)} = -|B|$, and $B^{(1)} = |B|$ and will not be considered further.
\end{remark}

It is not difficult to see from this
definition that skew Brownian paths inherit almost sure continuity from
that of Brownian motion $B$. Moreover,
let ${\cal F}_t:= \sigma\{|B(s)|\dvtx  0\le s\le t\}\vee\sigma\{
A_1,A_2,\ldots\}$, therefore, ${\cal F}_t \supseteq\sigma(B^{(\alpha
)}(s)\dvtx  s\le t)$.
For $0\le s < t$ and
a nonnegative,
measurable function $g$, one may use the
Markov property of $|B|$, and the independence of $|B|$
from the i.i.d. sign changes $A_1,A_2,\ldots\,$, to check that
\[
\EXP\bigl\{g\bigl(B^{(\alpha)}(t)\bigr) | {\cal F}_s\bigr\} =
\EXP\bigl\{g\bigl(B^{(\alpha)}(t)\bigr)| B^{(\alpha)}(s)\bigr\}.
\]
The Markov property of $B^{(\alpha)}$ follows since
${\cal F}_s\supseteq\sigma(B^{(\alpha)}(u)\dvtx  u\le s)$.

Using the strong Markov property for Brownian motion, \citet
{walsh1} calculated the transition
probabilities $p^{(\alpha)}(t,x,y)$ for $\alpha$-skew Brownian motion
as given by
%
%
\begin{eqnarray}
\label{eqskewtrans}
&&p^{(\alpha)}(t,x,y)\nonumber\\[-8pt]\\[-8pt]
&&\quad= \cases{\displaystyle\frac{1}{\sqrt{2\pi t}}
e^{-{(y-x)^2}/({2t})} \vspace*{2pt}\cr
\quad{}+ \displaystyle\frac{(2\alpha
-1)}{\sqrt{2\pi t}} e^{-{(y+x)^2}/({2t})}, \vspace*{2pt}\cr
\quad \hspace*{102pt}\mbox{if $x>0, y>0$},
\vspace*{2pt}\cr
\displaystyle\frac{1}{\sqrt{2\pi t}} e^{-{(y-x)^2}/({2t})} \vspace*{2pt}\cr
\quad{}- \displaystyle\frac{(2\alpha-1)}{\sqrt{2\pi t}}e^{-
{(y+x)^2}/({2t})}, \vspace*{2pt}\cr
\quad \hspace*{102pt}\mbox{if
$x<0, y<0$},
\vspace*{2pt}\cr
\displaystyle\frac{2\alpha}{\sqrt{2\pi t}} e^{-{(y-x)^2}/({2t})},
\quad\hspace*{12.3pt} \mbox{if $x\leq0, y>0$},
\vspace*{2pt}\cr
\displaystyle\frac{2(1-\alpha)}{\sqrt{2\pi t}} e^{-{(y-x)^2}/({2t})},
\quad \mbox{if $x\geq0, y<0$}.}\hspace*{-25pt}
\nonumber
\end{eqnarray}

Now,
since a Markov process is uniquely determined by its transition
probabilities and initial distribution, it is a simple
matter to use a change of variable transformation to check that the
physical skew diffusion $X^*$
is a particular (rescaling) function of
a skew Brownian motion with a particular transmission coefficient
$\alpha= \alpha^*$. In particular,
the above may be summarized as the following:
%
\begin{theorem}
\label{skewdiffderived}
Define the physical skew diffusion process $X^* = \{X^*(t)\dvtx  t\ge0\}$
by
%
%
\begin{eqnarray}
X^*(t) = s_{\sqrt{D}}\bigl(B^{(\alpha^*)}(t)\bigr),\nonumber\\[-8pt]\\[-8pt]
\eqntext{\displaystyle  t\ge0, \alpha^* =
\frac{\sqrt{D^+}}{\sqrt{D^+} + \sqrt{D^-}},}
\end{eqnarray}
where
%
%
\begin{equation}
\label{Dscalefctn} s_{\sqrt{D}}(x) = \cases{ \sqrt{D^+}x, & if $x > 0$,
\vspace*{2pt}\cr
\sqrt{D^-}x, & if $x \leq0$.}
\end{equation}
Then $X^*$ is the diffusion on $\R$ with transition probabilities given
by (\ref{eqskewdifftrans})
started at zero.
\end{theorem}

\begin{remark} As previously noted, the self-\break adjointness property that
results from the conservative interface
condition (\ref{conserveinterface}) may be viewed as a symmetry of the
transition probabilities
(\ref{eqskewdifftrans}) of the physical skew diffusion.
Although one sees by inspection that the transition probabilities (\ref
{eqskewtrans}) of skew Brownian motion
are \textit{not} symmetric in the sense of (\ref{eqskewdifftrans}), using the
strong Markov property of skew Brownian motion,\footnote{One may check
that for fixed
$y$, the transition probabilities of the skew Brownian motion are
continuous in the
backward variable~$x$, that is, a \textit{Feller property} holds. As a
consequence of this and the
sample path continuity, the \textit{strong Markov property} for skew
Brownian motion also
follows.} one has that
$r_\alpha(B^{(\alpha)})$ is a martingale, where
${r}_\alpha(x) = \alpha x{\mathbf1}_{(-\infty,0]}(x)
+ (1-\alpha)x{\mathbf1}_{[0,\infty)}(x), x\in\R$; see \citet{walsh1}.
\end{remark}

\begin{remark} A similar, albeit somewhat more technical, procedure may
be developed for
%
%
\begin{eqnarray}
\label{1adskewpde} \pderiv{u} {t}(t,x) &=& \frac{1}{2}D(x)\,\pderiv[2]{u}
{x}(t,x) \nonumber\\
&&{}- v \,\pderiv {u} {x}(t,x)\quad (x\neq0),
\nonumber\\[-8pt]\\[-8pt]
\quad D^+\,\pderiv{u} {x}\bigl(t,0^+\bigr) &=& D^-\,\pderiv{u} {x}\bigl(t,0^-
\bigr),\nonumber\\
\eqntext{u\bigl(t,0^+\bigr) = u\bigl(t,0^-\bigr), t > 0,}
\end{eqnarray}
for a constant drift $v$ by a change of measure that converts the
problem into one of an
elastic skew Brownian motion; see \citet{AppuhamilageAAP} for details.
\end{remark}

While we have perhaps traced the most direct route from the p.d.e. model
(\ref{1dskewpde}) to probability theory, several others naturally
emerge. Some of these are summarized in the next section.

\section{Alternative Mathematical Descriptions of Physical Skew
Diffusion}\label{sec3} \label{SectionAlternatives}

In this section four equivalent approaches to represent the Markov
process associated with (\ref{1dskewpde}) are provided as alternatives
to the construction in terms of excursions of reflected Brownian motion
paths. Each of these provides additional mathematical tools in which to
gainfully address diverse problems involving dispersion in the presence
of the conservative interface condition (\ref{conserveinterface}). We
begin with perhaps the most mathematically technical framework, that of
Dirichlet forms, which can be useful for existence theory and for weak
formulations used in developing numerical methods. This subsection can
certainly be skimmed on first reading. The overall section progresses
to the least technical framework of skew random walks and is followed
by a subsection addressing a more general class of (nonconservative)
interface conditions that will be seen to arise naturally in certain
physical, biological and ecological dispersion contexts.

\subsection*{Dirichlet Forms}

Below we outline the procedure leading to a semigroup framework for
$X^*$ via Dirichlet forms theory, and refer the reader to the more
comprehensive references by \citet{Fukushima94}, \citet
{marockner} or the recent \citet{chenfuki}.

To set up the analytical framework, let $u$ be a solution to problem
(\ref{1dskewpde}) and consider the following variational form of the
evolution equation in $L^2(\R)$:
%
%
\begin{eqnarray}
&&
\pderiv{} {t} \int_\R u(t,x) v(x) \sud x \nonumber\\
&&\quad= - \int
_\R\frac{1}{2} D(x) \,\pderiv{u} {x}(t,x) \,\pderiv{v}
{x}(x) \sud x,\\
\eqntext{v \in C_c^\infty(\R).}
\end{eqnarray}
The associated process is obtained by identifying a semigroup generated
by the bilinear form
%
%
\begin{eqnarray}
\label{defE} \E(u,v) = \int_\R\frac{1}{2} D(x)\,
\pderiv{u} {x}(x) \,\pderiv {v} {x}(x) \sud x,\nonumber\\[-8pt]\\[-8pt]
\eqntext{u,v \in C_c^\infty(
\R),}
\end{eqnarray}
in some Hilbert space. For the case of $D$ given by (\ref{1dskewpde}),
standard considerations show that $\E$ is ``closable'' on $L^2(\R)$,
namely, it extends to a closed bilinear form (also denoted by $\E$) on
$L^2(\R)$ with domain $\Dom(\E) = H^1(\R)$, the Sobolev space of $L^2$
functions whose generalized derivatives are also square integrable
functions. Here, ``closed'' means that $\E(u,u) \geq0$ for all $u \in
\Dom(\E)$ and that $\Dom(\E) = H^1(\R)$ is a Hilbert space with the
inner product $\E_1(u,v):= \E(u,v) + (u,v)_{L_2(\R)}$. The bilinear
form $(\E,\Dom(\E))$ is ``coercive'' since $\E(u,v)^2 \leq\E_1(u,u) \E
_1(v,v)$ and a ``Dirichlet form'' since
%
%
\begin{eqnarray}
\label{DirichletProp} \E(u,u) \leq\E\bigl(u^+\wedge1,u^+\wedge1\bigr)\nonumber\\[-8pt]\\[-8pt]
\eqntext{\mbox{for all } u,v \in \Dom(\E).}
\end{eqnarray}
Finally, $(\E,\Dom(\E))$ is ``regular'' since $\Dom(\E) \cap C_c(\R)$ is
dense in $\Dom(\E)$ with respect to the norm $u \mapsto\E(u,u)^{1/2}$.
For such a form, there exists a unique closed, negative definite,
linear operator $(\Ab,\Dom(\Ab))$ that satisfies the resolvent
conditions of the Hille--Yosida theorem for generating the
appropriate semigroup; namely, $(\lambda- \Ab)(\Dom(\Ab)) = L^2(\R)$,
$\lambda> 0$. This operator is given by
%
%
\begin{eqnarray}
\Dom(\Ab) &\subset&\Dom(\E) \quad\mbox{and}\nonumber\\
\E(u,v) &=& (-\Ab u,v) \\
\eqntext{\mbox{for all } u \in
\Dom(\Ab), v \in\Dom(\E).}
\end{eqnarray}
Integration by parts on (\ref{defE}) yields
%
%
\begin{eqnarray}
\label{AL2} &&\Ab f = \frac{1}{2}D \,\pderiv[2]{f} {x},
\nonumber\\
&&\Dom(\Ab) \nonumber\\[-8pt]\\[-8pt]
&&\quad= \biggl\{ f \in H^1(\R) \cap H^2
\bigl(\R^+\bigr) \cap H^2\bigl(\R^-\bigr)\dvtx\nonumber\\
&&\hspace*{57.5pt}  D^+ \,\pderiv{f} {x}
\bigl(0^+\bigr) = D^- \,\pderiv{f} {x}\bigl(0^-\bigr) \biggr\},\nonumber
\end{eqnarray}
where $H^2(\R^\pm)$ denote the respective
Sobolev spaces for twice (generalized) differentiable functions on $\R
^\pm$.
The operator $(\Ab, \Dom(\Ab))$ is the infinitesimal generator of a
strongly continuous contraction semigroup $\{T_t\dvtx t\ge0\}$
on $L^2(\R)$ which is also sub-Mar\-kovian since $(\E, \Dom(\E))$
satisfies the Dirichlet form property (\ref{DirichletProp}). Note also
that the conservative interface condition is encoded in
$\Dom(\Ab)$ and makes the operator $\Ab$ self-adjoint.
The family of transition probabilities
$p^*(t,\cdot, \cdot)$, $t>0$ are recovered from the semigroup via, for bounded
$A\in{\cal B}({\mathbb R})$,
%
%
\begin{equation}
p^*(t,x,A) = T_t\IND{A}(x).
\end{equation}

Since for $f \in\Dom(\Ab)$, the unique solution in\break  $C([0,\infty)) \cap
\Dom(\Ab)$ to $\pderiv{u}{t} = \Ab u$, $u(0,x) = f(x)$ is $u(t,x) = \int_\R p^*(t,x,\ud y) f(y)$ with $p^*$ given in (\ref{eqskewdifftrans}),
then $T_t f = \int_\R p^*(t,\cdot,y) f(y) \sud y$ almost everywhere for
any $f \in L^2(\R)\cap L^\infty(\R)$. The process
$X^*$ therefore has transition probabilities given by (\ref{eqskewdifftrans}).




\subsection*{Feller's Natural Scale and L\'evy's Time Change of
Brownian Motion}

In this subsection we outline the procedure leading to the
characterization of $X^*$ via Feller's natural scale and L\'evy's time
change of Brownian motion, and refer the reader to the more general
references \citet{revuzyor},
\citet{karshrevebook}, \citet{bhattway} and others.

Let $a<x<b$ be arbitrary. The scale measure $s_*$ of $X^*$ is the
unique (up to a multiplicative constant) measure that satisfies
%
%
\begin{equation}
\label{Ptaus} \P_x\bigl(\tau^*_{(a,b)} = H^*_b
\bigr) = \frac{s_*((a,x))}{s_*((a,b))},
\end{equation}
where $H^*_b:= \inf\{t\ge0\dvtx  X^*(t) = b\}$ denotes the \textit{hitting
time} of $b$,
and $\tau^*_{(a,b)}:= \inf\{t\ge0\dvtx  X^*(t)\in\{a,b\}\}$ is the \textit{escape time} from
the interval $(a,b)$ for $X^*(0)\in(a,b)$.
The speed measure $m_*$ of $X^*$ is the unique Radon measure on Borel
subsets of $\R$ such that
%
%
\begin{equation}
\label{EXPtaum} \EXP_x \tau^*_{(a,b)} = \int
_a^b G_{a,b}(x,y)m_*(\ud y),
\end{equation}
where the so-called Green's function $G$ of $X^*$ is given by
%
%
\begin{equation}\quad
G_{a,b}(x,y) = \frac{s_*((x \wedge y,a)) s_*((b,x \vee y))}{s_*((a,b))}.
\end{equation}
%

The process $X^*$ has speed and scale measures with piecewise constant
density. Let $m_*(\ud x) = m_*'(x) \sud x$, $s_*(\ud x) = s_*'(x) \sud
x$ with
%
%
\begin{eqnarray}\qquad
m_*'(x) &=& m_*^- \IND{(-\infty,0)}(x) + m_*^+
\IND{(0,\infty)}(x),\nonumber\\[-8pt]\\[-8pt]
s_*'(x) &=& s_*^- \IND{(-\infty,0)}(x) + s_*^+
\IND{(0,\infty)}(x).\nonumber
\end{eqnarray}

To determine the constants, let $(\A,\Dom(\A))$ denote the restriction
of the operator $(\Ab,\Dom(\Ab))$ in (\ref{AL2}) to the space of
bounded continuous functions:
%
%
\begin{eqnarray}
\label{ACb}\quad &&\A f = \frac{1}{2}D \,\pderiv[2]{f} {x},
\nonumber\\
&&\Dom(\A) \nonumber\\[-8pt]\\[-8pt]
&&\quad= \biggl\{ f \in C_b(\R) \cap C^2
\bigl(\R^+\bigr) \cap C^2\bigl(\R^-\bigr)\dvtx\nonumber\\
&&\hspace*{54.5pt}  D^+ \,\pderiv{f} {x}
\bigl(0^+\bigr) = D^- \,\pderiv{f} {x}\bigl(0^-\bigr) \biggr\}.\nonumber
\end{eqnarray}
Then $(\A,\Dom(\A))$ is also given by
%
%
\begin{eqnarray}
\label{Admds} \A f &=& \deriv{} {m_*} \,\deriv{} {s_*}f,\nonumber\\[-8pt]\\[-8pt]
\hspace*{28pt}\Dom(\A) &=& \biggl\{f \in
C_b(\R)\dvtx  \deriv{f} {s_*} \mbox{ is continuous} \biggr\},\nonumber
\end{eqnarray}
where $\deriv{f}{s_*} = \deriv{f}{x}\,\deriv{x}{s_*} = \deriv{f}{x}/\deriv
{s_*}{x}$.
Matching the expressions for $\A$ above, one arrives at
%
%
\begin{equation}\qquad
s_*^+ = \frac{c}{D^+},\quad s_*^- =\frac{c}{D^-},\quad m_*^+ = m_*^- =
\frac{2}{c},
\end{equation}
where $c$ is any positive constant which we set equal to one for convenience.

Within this framework (see \citecs{revuzyor}, page 310) one has the following:
%
\begin{theorem}
Physical skew diffusion $X^*$ with $D$ given by (\ref{1dskewpde}) is
the unique Feller process on $\R$ with speed and scale measures
$m_*(\ud x) = m_*'(x) \sud x$, $s_*(\ud x) = s_*'(x) \sud x$ with
densities given by
%
%
\begin{eqnarray}\qquad
m_*'(x) &=& 2,\nonumber\\[-8pt]\\[-8pt]
s_*'(x) &=& \frac{1}{D^-} \IND{(-
\infty,0)}(x) +\frac{1}{D^+} \IND {(0,\infty)}(x).\nonumber
\end{eqnarray}
\end{theorem}

To see the propagation of local interface effects on the global
features within this framework, let $a>0$ and $0<\eps<a$, and use (\ref
{Ptaus}), (\ref{EXPtaum}) to obtain
%
%
\begin{eqnarray}
\EXP_a \tau^*_{(a-\eps,a+\eps)} &=&
\frac{\eps^2}{D^+},\nonumber\\[-8pt]\\[-8pt]
\P_a\bigl(\tau ^*_{(a-\eps,a+\eps)} &=& H^*_{a-\eps}\bigr) =
\frac{1}{2}\nonumber
\end{eqnarray}
as expected, since starting at $a$, $X^*$ must behave like a diffusion
process with diffusion coefficient $D^+$ up to the hitting time $H_0 >
\tau_{(a-\eps, a+\eps)}$. On the other hand, for the process starting
at the interface at $x=0$, the effects of the heterogeneity are
depicted~by
%
%
\begin{eqnarray}
\label{skewnessX} \EXP_0 \tau^*_{(-\eps,\eps)} &=& \frac{2\eps^2}{D^+ +
D^-},\nonumber\\[-8pt]\\[-8pt]
\P _0\bigl(\tau^*_{(-\eps,\eps)} &=& H^*_{-\eps}\bigr) =
\frac{D^-}{D^+ + D^-}.\nonumber
\end{eqnarray}
Namely, the interface $x=0$ ``skews'' the process, making it more likely
to exit the symmetric interval $(-\eps,\eps)$ through the endpoint with
highest diffusion coefficient value.

For a path-wise representation of the process $X^*$ one may proceed by
L\'evy's time change of Brownian motion as follows. Let $B = \{B(t)\dvtx
{t\geq0}\}$ denote
canonical standard Brownian motion on
$(\Omega, \mathcal{F},\break  \{\P_x\}_{x \in\R})$ and consider the following
additive functional
%
%
\begin{equation}
\phi(r) = \int_0^r \frac{1}{2}
\frac{m_*'(B(t))}{s_*'(B(t))} \sud t,\quad r \geq0.
\end{equation}
Let $T$ be the inverse of $\phi$, $T(t) = \inf\{s \geq0\dvtx  \phi(s) = t
\}$, then the process $X^*$ has the following representation:
%
%
\begin{equation}
\label{skdiffRTC} X^*(t) = s_*^{-1}\bigl(B\bigl(T(t)\bigr)\bigr),\quad t
\geq0,
\end{equation}
where $s_*^{-1} $ denotes the inverse of the function $x \mapsto
s_*((0,x))$, $s_*^{-1}(x) = D(x)x$.


The representation obtained in (\ref{skdiffRTC}) can be simplified
further in order to write $X^*$ as a function of a continuous
martingale. Consider the time-change
$Y(t) = B(T(t))$, where
%
%
\begin{eqnarray}
T(t) &=& \int_0^{T(t)} 2 \frac{s_*'(B(r))}{m_*'(B(r))} \sud
\phi(r) \nonumber\\[-8pt]\\[-8pt]
&=& \int_0^t 2\frac{s_*'(Y(\rho))}{m_*'(Y(\rho))} \sud
\rho\nonumber
\end{eqnarray}
and, therefore, the quadratic variation of $Y$ is\break  $\langle Y \rangle(t) =
T(t)$. But since $\phi(r)$ is continuous, increasing and finite, then
so is $T$. Therefore (see Karatzas and Shreve, \citeyear{karshrevebook}, Theorem 4.2), there
exists a probability space\vadjust{\goodbreak} $(\tilde\Omega, \tilde{\mathcal{F}},
\{\tilde\P_x\}_{x \in\R})$ extending $(\Omega, \mathcal{F},\break
\{\P_x\}_{x \in\R})$ and with a Brownian motion $\tilde B$ defined
such that
%
%
\begin{eqnarray}
\label{YtTt} Y(t) &=& \int_0^t Z(r) \sud\tilde
B(r),\nonumber\\[-8pt]\\[-8pt]
T(t)&=&\langle Y\rangle(t) = \int_0^t
Z^2(r) \sud r,\quad \tilde\P\mbox{-a.s.}\nonumber
\end{eqnarray}
for some measurable adapted process $Z$.
It follows from (\ref{YtTt}) that $Z(t) =
[2s_*'(Y(t))/m_*'(Y(t))]^{1/2}$. So one arrives at the following
representation in terms of the martingale $Y$:
%
%
\begin{eqnarray}
\label{skdiffMart} X^*(t) &=& D\bigl(Y(t)\bigr) Y(t),\nonumber\\[-8pt]\\[-8pt]
Y(t) &=& \int
_0^t \frac{1}{\sqrt {D(Y(r))}} \sud\tilde{B}(r).\nonumber
\end{eqnarray}

\subsection*{Stochastic Calculus and Local Time}

The representation of $X^*$ in (\ref{skdiffMart}) again makes it
evident that whenever the process is away from the interface, the
trajectories of $X^*$ can be obtained by simply re-scaling those of
Brownian motion by the square root of the appropriate diffusion
coefficient. However, it does not reveal the behavior of $X^*$ at $x=0$
and, in particular, the skewness property (\ref{skewnessX}). This
property must be produced by the effect of the jump in the value of the
diffusion coefficient over the trajectories, during the ``time'' a
particle occupies the interface. In order to quantify this effect then,
one is naturally led to consider the properties of the local time of $X^*$.

We now briefly give some necessary background on the theory of local
time for continuous semi-martingales. The reader is referred to
\citet{revuzyor} for the general theory followed here. Given a
continuous semimartingale $X(t) = M(t) + V(t)$, where $M$ is a
martingale and $V$ is an increasing process, we define its local time
process $\Lp{X}$ via
%
%
\begin{eqnarray}
\label{defL} \bigl|X(t) - a\bigr| &=& \bigl|X(0)-a\bigr| \nonumber\\
&&{}+ \int_0^t
\sign_-\bigl(X(s)-a\bigr) \sud X(s)\\
&&{} + \Lp{X}(t,a)\nonumber
\end{eqnarray}
with the convention $\sign_- = \IND{(0,\infty)} - \IND{(-\infty,0]}$.
What we are calling local time in this paper is sometimes referred to
in the literature as
\textit{right local time} since it satisfies almost surely
%
%
\begin{equation}
\label{altDef}\hspace*{28pt} \Lp{X}(t,a) = \lim_{\varepsilon\downarrow0} \frac{1}{\varepsilon} \int
_0^t \IND{[a,a+\varepsilon)}\bigl(X(r)\bigr) \sud
\langle X \rangle(r).
\end{equation}
The function $(t,a) \mapsto\Lp{X}(t,a)$ can be taken continuous in $t$
and cadlag in $a$ and its jumps\vadjust{\goodbreak} are given by
%
%
\begin{eqnarray}
\label{jumpL} &&
\Lp{X}(t,a) - \Lp{X}\bigl(t,a^-\bigr) \nonumber\\
&&\quad= 2 \int_0^t
\IND{\{a\}}\bigl(X(s)\bigr) \sud X(s) \\
&&\quad= 2 \int_0^t
\IND{\{a\}}\bigl(X(s)\bigr) \sud V(s).\nonumber
\end{eqnarray}
In particular, if $X= M$ is a local martingale, then $(t,a) \mapsto\Lp
{X}(t,a)$ can be taken bi-continuous.

The following basic formulae encompass the most significant properties
of local time:
\begin{longlist}
\item[\textit{It\^o--Tanaka formula}.] If $f$ is a difference of convex functions
and $f'_-$ denotes its left derivative, then
%
%
\begin{eqnarray}
\label{IItoTanaka}\hspace*{28pt} f\bigl(X(t)\bigr) &=& f\bigl(X(0)\bigr) + \int
_0^t f'_-\bigl(X(s)\bigr) \sud
X(s)\nonumber\\[-8pt]\\[-8pt]
&&{} + \frac{1}{2} \int_\R \Lp{X}(t,x)
f''(\ud x).\nonumber
\end{eqnarray}
\item[\textit{Occupation times formula}.] For any positive Bo-\break rel---measurable
function $F$,
%
%
\begin{equation}
\label{occupTimes}\hspace*{28pt} \int_0^t F\bigl(X(s)\bigr)
\sud\langle X\rangle(s) = \int_\R F(x) \Lp{X}(t,x) \sud x.
\end{equation}
\item[\textit{Left-side local time}.] If $X$ is a continuous semimar\-tingale, then
almost surely
%
%
\begin{equation}
\label{altDefLeft} \Lp{X}\bigl(t,a^-\bigr) = \lim_{\varepsilon\downarrow0}
\frac{1}{\varepsilon} \int_0^t \IND{(a-
\varepsilon,a)}\bigl(X(r)\bigr) \sud\langle X \rangle(r).\hspace*{-34pt}
\end{equation}
\end{longlist}

As a first step consider the representation of $X^*$ in (\ref
{skdiffMart}) as a nonsmooth function of the martingale~$Y$. Applying
the It\^o--Tanaka formula to the function $f(x) = s_*^{-1}(x) = xD(x)$
and using the representation (\ref{skdiffMart}) of $Y$ in terms of a
Brownian motion $B$, one gets
\begin{eqnarray*}
X^*(t) &=& \int_0^t D\bigl(Y(r)\bigr) \sud
Y(r) + \frac{1}{2}\bigl(D^+ - D^-\bigr)L^Y(t,0)
\\
&=& \int_0^t \sqrt{D\bigl(X^*(r)\bigr)} \sud
B(r) \\
&&{}+ \frac{1}{2}\bigl(D^+ - D^-\bigr)L^Y(t,0).
\end{eqnarray*}
The local time of $Y$ can be related to the local time $\Lp{*}$ of
$X^*$ using
(\ref{altDef}). Note first that $\langle X^*\rangle(r) = D^2(X(r))
\langle Y\rangle(r) $ for $r\geq0$, then write
%
%
\begin{eqnarray}\quad
&&\Lp{Y}(t,a) \nonumber\\
&&\quad= \lim_{\varepsilon\downarrow0} \frac{1}{\varepsilon} \int
_0^t \IND{[0,\varepsilon)}\bigl(Y(r)\bigr) \sud
\langle Y \rangle(r)
\nonumber\\[-8pt]\\[-8pt]
&&\quad= \lim_{\varepsilon\downarrow0} \frac{ 1}{(D^+)^2\varepsilon} \int_0^t
\IND{[0,D^+ \varepsilon)}\bigl(X(r)\bigr)\sud\langle X \rangle(r) \nonumber\\
&&\quad=
\frac{1}{D^+} \Lp{*}(t,0).\nonumber
\end{eqnarray}
We have arrived at the following representation of~$X^*$.

\begin{theorem}\label{thmStrongSol}
For a given $X(0)$, and on any filtered probability space carrying a
Brownian motion $B$, physical skew diffusion $X^*$ is the unique strong
solution to the following stochastic differential equation:
%
%
\begin{eqnarray}
\label{skewDiffSDE} X(t) &=& X(0) + \int_0^t
\sqrt{D\bigl(X(r)\bigr)} \sud B(r) \nonumber\\[-8pt]\\[-8pt]
&&{}+ \frac{(D^+ -
D^-)}{2D^+}L^X(t,0).\nonumber
\end{eqnarray}
\end{theorem}

Equation (\ref{skewDiffSDE}) is a stochastic differential equation in
terms of the local time of the unknown process. \citet{LeGall}
studied the problem of existence and uniqueness of solutions for
equations of this type and proved Theorem \ref{thmStrongSol}. In fact,
he considered a larger
set of equations which we review here for its relevance with regard to
more general solute transport problems. For the sake of consistency
with other parts of the present paper, we
summarize his analysis in terms of the right local time defined in (\ref
{defL}), in place of the symmetric local time used in \citet{LeGall}.

Consider a finite signed measure $\nu$
such that\break  $\nu(\{x\}) < 1$ for all $x \in\R$, and let $\nu^c$ be its
continuous part. Also, let $\vp$ be a right-continuous function of
bounded variation that is also strictly positive and bounded away from
zero. \citet{LeGall} considered the following equation:
%
%
\begin{eqnarray}
\label{LeGallSDE} X(t) &=& X(0) + \int_0^t \vp
\bigl(X(r)\bigr) \sud B(r) \nonumber\\[-8pt]\\[-8pt]
&&{}+ \frac{1}{2} \int_\R
L^X(t,x) \nu(\ud x).\nonumber
\end{eqnarray}
In the case of equation (\ref{skewDiffSDE}), $\vp= \sqrt{D}$ and $\nu
= \frac{D^+-D^-}{D^+} \delta_0$, in particular,\vspace*{1pt} $\nu^c \equiv0$. The
key to the analysis is to relate equation (\ref{LeGallSDE}) to a
stochastic differential equation without a local time term. In fact,
\citet{LeGall} shows that if $f_\nu$ is a right-continuous
function satisfying
%
%
\begin{equation}\quad
f_\nu\bigl(x^-\bigr) \nu(\ud x) + f_\nu'(
\ud x) = 0,\quad x \in\R,
\end{equation}
and $F_\nu(x) = \int_{-\infty}^x f_\nu(y) \sud y$, then a process $X$
is a solution to (\ref{LeGallSDE}) if and only
%
%
\begin{eqnarray}
Y(t) &=& F_\nu\bigl(X(t)\bigr) \quad\mbox{satisfies}\nonumber\\[-8pt]\\[-8pt]
Y(t) &=& \int
_0^t f_\nu \bigl(X(r)\bigr) \vp
\bigl(X(r)\bigr) \sud B(r).\nonumber
\end{eqnarray}
Moreover, it is easy to show that the function $f_\nu$ is given by
%
%
\begin{equation}
f_\nu(x) = \exp\bigl(-\nu^c\bigl((-\infty,x]\bigr)\bigr) \prod
_{y \leq x}\bigl(1-\nu\bigl(\{y\}\bigr)\bigr).\hspace*{-32pt}
\end{equation}
Not surprisingly, when this procedure is applied to equation (\ref
{skewDiffSDE}), one gets $f_\nu= 1/D$ and recovers problem (\ref
{skdiffMart}). The existence and uniqueness of strong solutions to (\ref
{skewDiffSDE}) follows then from the corresponding result for (\ref
{skdiffMart}) which was established in \citet{Nakao}, and recently
generalized substantially in\break  \citet{Prokaj} and \citet
{Fernholz2012fk}.

It is important to note that the representation (\ref{skewDiffSDE})
gives the decomposition of the continuous semi-martingale $X^*$ as the
sum of a continuous local martingale and an increasing process. It
follows that the local time of $X^*$ is not continuous at $x=0$. In
fact, by (\ref{jumpL}), we can compute
%
%
\begin{equation}
\frac{L^*(t,0)}{L^*(t,0^-)} = \frac{D^+}{D^-},\quad t \geq0.
\end{equation}
This, however, cannot be interpreted as skew diffusion ``spending more
time'' on either side of the interface. To see this, we can use the
alternative definitions (\ref{altDef}), (\ref{altDefLeft}) of right and
left local time. For $X^*$, (\ref{skewDiffSDE}) gives the quadratic
variation $\langle X^*\rangle(t) = \int_0^t D(X(r)) \sud r$. In
particular, the ratio between the time a particle spends just above the
interface $x=0$ up to time $t$ and the time it spends just below that
interface is
%
%
\begin{eqnarray}
&&
\lim_{\eps\downarrow0} \frac{\int_0^t \IND{[0,\eps)}(X^*(r)) \sud r} {
\int_0^t \IND{(-\eps,0]}(X^*(r)) \sud r}\nonumber\\
\hspace*{20pt}&&\quad = \lim_{\eps\downarrow0}
\frac{(\int_0^t \IND{[0,\eps )}(X^*(r)) \sud\langle X^*
\rangle(r))/{D^+}}{(\int_0^t \IND{(-\eps,0]}(X^*(r)) \sud\langle
X^* \rangle(r))/{D^-}} = 1,\\
\eqntext{t\geq0.}
\end{eqnarray}
We will return to such matters in the context of applications.

\subsection*{Discrete and Numerical Approximations}

\textit{Skew random walk}
is a natural discretization of skew Brownian motion
defined as follows.
%
\begin{definition}\label{skewrw}
The \textit{$\alpha$-skew random walk} is a discrete Markov chain $\{Y_n\dvtx
n=0,1,2,\ldots\}$ on the integers $\mathbb{Z}$
having transition probabilities
\[
p^{(\alpha)}_{ij} = \cases{\frac{1}{2}, & if $i \neq0, j=i
\pm1$,
\vspace*{2pt}\cr
\alpha, & if $i = 0, j=1$,
\vspace*{2pt}\cr
1-\alpha, & if $i = 0, j= -1$.}
\]
\end{definition}

Convergence of the distribution at a fixed time point was first
announced in \citet{harrisonshepp}, where they indicated a
``fourth moment\break  proof'' along the lines of that given for a simple
symmetric random walk (i.e., $\alpha= 1/2$) based on convergence of
finite-dimensional distributions. However, proving tightness is quite
laborious and tricky due to the lack of independence of the increments.
A~full proof was given in \citet{brookschacone}. The remainder of
this section describes an approach based on the Skorokhod embedding
method within this more specialized framework. A more general
functional central limit theorem is given in
\citet{chernyshiryor}.

\begin{lemma}[(Discrete excursion representation)]
\label{discreteexcursionrep}
Let ${S} = \{S_n\dvtx  n =0,1,\ldots\}$ be a simple symmetric random
walk starting at $0$, and let $\tilde{J}_{\pi_1},\tilde{J}_{\pi_2},\ldots$
denote an enumeration of the excursions of ${S}$
away from zero for a fixed but arbitrary permutation
$\pi$
of the natural numbers. In particular, $|S_n| > 0$ if $n\in\tilde
{J}_{\pi_k}$.
Define
\[
S^{(\alpha)}_0 = 0,\quad S_n^{(\alpha)} = \sum
_{k=1}^\infty\IND{\tilde{J}_{\pi_k}} (n)
\tilde{A}_k|S_n|,\quad n \ge1,
\]
where $\tilde{A}_1,\tilde{A}_2,\ldots$ is an i.i.d. sequence of
Bernoulli $\pm1$-random variables, independent of ${S}$, with
$\P(\tilde{A}_1 = 1) = \alpha$.
Then ${S}^{(\alpha)}$ is distributed as an $\alpha$-skew random
walk.
\end{lemma}

Define the polygonal random function $S^{(\alpha,n)}$ on $[0,1]$ as follows:
%
%
\begin{eqnarray}
\qquad S^{(\alpha,n)}(t):=\frac{S^{(\alpha)}_{k-1}}{\sqrt{n}} - \frac
{S_k^{(\alpha)}-S_{k-1}^{(\alpha)}}{\sqrt{n}} \biggl( t -
\frac{k-1}{n} \biggr),\nonumber\\[-8pt]\\[-8pt]
\eqntext{\displaystyle  t \in \biggl[\frac{k-1}{n},\frac{k}{n}
\biggr], 1 \leq k \leq n.}
\end{eqnarray}

That is, $S^{(\alpha,n)}(t) = \frac{S^{(\alpha)}_k}{\sqrt{n}}$ at points
$t= \frac{k}{n}$ $(0\le k\le n)$, and
$t\mapsto S^{(\alpha,n)}(t)$ is\vspace*{1pt} linearly interpolated
between the endpoints of each interval
$ [\frac{k-1}{ n},\frac{k}{n} ]$.

Let us recall that by an application of the
Skorokhod embedding theorem (e.g., see Bhattacharya and Waymire, \citeyear{bcpt}),
there is a sequence of
times $T_1 < T_2 <\cdots$ and a Brownian motion $\{B(t)\dvtx  t \geq0\}$
such that
$B(T_1)$ has a symmetric Bernoulli $\pm1$-distribution, and
$B(T_{i+1}) - B(T_i)$ $(i\ge0)$ are i.i.d.
with a symmetric $\pm1$-distribution. Moreover,
$T_{i+1} - T_i$ $(i\ge0)$ are
i.i.d. with mean one. With this one may check the following:
%
\begin{lemma}
The discrete parameter stochastic process
$\tilde S^{(\alpha)}_0 = 0$,
$\tilde S^{(\alpha)} _m:= B^{(\alpha)}(T_m)$, $m = 1,2,\ldots\,$,
is distributed as an $\alpha$-skew random walk.
\end{lemma}
Now it is a rather straightforward exercise to prove the following
theorem as an application of the Skorokhod embedding theorem, similar
to that for weak convergence of the simple random walk to Brownian
motion found in \citet{bcpt} and many other references.
%
\begin{theorem}
\label{thw5111}
$S^{(\alpha,n)}$ converges in distribution to the $\alpha$-skew Brownian
motion $B^{(\alpha)}$
as $n\to\infty$.
\end{theorem}
Since the rescaling function is continuous,
it follows that the rescaled skew random walks converge in distribution
to the physical skew diffusion. That is, recalling the definition of
$s_{\sqrt{D}}$
at (\ref{Dscalefctn}), one has the following:
%
\begin{cor}
The (polygonal) random walks $X_n^*, n\ge1$,
defined by
\[
X_n^*(t) = s_{\sqrt{D}}\bigl(\tilde S^{(\alpha^*,n)}(t)
\bigr),\quad t\ge0, n = 1,2,\ldots,
\]
converge weakly to the physical skew diffusion
process $X^*$ on $C[0,\infty)$.
\end{cor}

The convergence of the discretized process opens the door to numerical
simulation schemes. Two important alternatives to numerical methods are
naturally suggested, namely, numerical solutions to the p.d.e. (\ref{1dskewpde})
and/or numerical solutions to the stochastic equation (\ref{skewDiffSDE}). The
self-adjoint character of the conservative interface conditions singles
out the numerical treatment of the equations in each of its
formulations (\ref{1dskewpde})
or (\ref{skewDiffSDE}). For example, standard off-the-shelf finite
difference methods provide numerical solutions to (\ref{1dskewpde}) for
the conservative interface condition. Similarly, in spite of the
presence of the local time term in (\ref{skewDiffSDE}), an
Euler/Muruyama method was designed by \citet{martineztalay} that
preserves the order of convergence of the Euler method when the
coefficients of the s.d.e. are smooth. They exploit the fact that in the
case of the conservative interface condition there is a one-to-one
piecewise linear transformation of the process that, with the aid of
the It\^o--Tanaka lemma, eliminates the local time term. As will be
emphasized in the next section and in subsequent examples to follow,
the conservative interface condition is only one of infinitely many
other possibilities of interest to applications that require new
approaches to numerical simulations of both the p.d.e. and s.d.e. Recently,
in \citet{lejay2012simulating}, \citet{Etore} and
in \citet{bokil2}, new
numerical methods, including both stochastic and deterministic
schemes, are developed
that apply to these more general interface conditions and restore the
order of convergence previously available for the more restrictive case
of the conservative interface.

\subsection*{General Interface Conditions}

The particular form of the interface condition (\ref
{conserveinterface}) arises naturally in the case of solute transport
as continuity of flux is imposed. However, as will be seen for
applications outside solute transport, the following more general
problem is also of interest:
%
%
\begin{eqnarray}
\label{1dLambdapde} \pderiv{u} {t} &=& \frac{1}{2} D \,\pderiv[2]{u}
{x},\nonumber\\
\lambda\,\pderiv{u} {x}\bigl(t,0^+\bigr) &=& (1-\lambda) \,\pderiv{u} {x}\bigl(t,0^-
\bigr),\\
\eqntext{u\bigl(t,0^+\bigr) = u\bigl(t,0^-\bigr), t >0,}
\end{eqnarray}
for some $0 < \lambda< 1$. The Markov process associated with problem
(\ref{1dLambdapde}) can be found using any of the techniques described
in this section (see Appuhamillage et~al.,
\citeyear{AppuhamAAPCorrect,AppuhamilageAAP}). In fact, skew Brownian
motion plays an important role here as provided by the following
extension of Theorem \ref{skewdiffderived}

\begin{theorem}
The Markov process associated with problem (\ref{1dLambdapde}) is
%
%
\begin{eqnarray}
\label{defXLambda} X(t) = s_{\sqrt{D}}\bigl(B^{(\alpha)}(t)\bigr),\qquad\qquad\nonumber\\[-8pt]\\[-8pt]
\eqntext{\displaystyle t\ge0, \alpha= \alpha(\lambda) = \frac{\lambda\sqrt{D^-}}{\lambda\sqrt{D^-}
+ (1-\lambda) \sqrt{D^+}},}
\end{eqnarray}
where $s_{\sqrt{D}}$ is given in (\ref{Dscalefctn}).
\end{theorem}

\begin{definition}
We refer to the Markov process associated to problem (\ref
{1dLambdapde}) as skew diffusion. In the special case
of the conservative interface condition for (\ref{1dskewpde}), we refer
to $X\equiv X^*$ as
the physical skew diffusion.
\end{definition}

\citet{Ouknine} characterizes skew diffusion processes as
solutions to a particular family of stochastic differential equations
of the form (\ref{LeGallSDE}). In particular, applying Tanaka's formula
gives that the process $X = s_{\sqrt{D}}(B^{(\alpha)})$ is a strong
solution to
%
%
\begin{eqnarray}\qquad
X(t) &=& X(0) + \int_0^t \sqrt{D\bigl(X(r)
\bigr)} \sud B(r) \nonumber\\[-8pt]\\[-8pt]
&&{}+ \frac{\alpha\sqrt {D^+} - (1-\alpha) \sqrt{D^-}}{2 \alpha\sqrt{D^+}}
L^X(t,0).\nonumber
\end{eqnarray}
The next theorem follows:

\begin{theorem}
Let $\gamma<1$, then the strong solution to
%
%
\begin{eqnarray}
X(t) &=& X(0) + \int_0^t \sqrt{D\bigl(X(r)
\bigr)} \sud B(r) \nonumber\\[-8pt]\\[-8pt]
&&{}+ \frac{\gamma}{2} L^X(t,0)\nonumber
\end{eqnarray}
is given by $X = s_{\sqrt{D}}(B^{(\alpha)})$ with
%
%
\begin{equation}
\label{aofgamma} \alpha= \frac{\sqrt{D^-}}{\sqrt{D^-} + \sqrt{D^+}(1-\gamma)}.
\end{equation}
\end{theorem}

Note that matching the formulae for $\alpha$ in (\ref{defXLambda}) and
(\ref{aofgamma}) gives
%
%
\begin{equation}
\lambda=\frac{1 }{ (2-\gamma)} \in(0,1),
\end{equation}
which expresses the discontinuities at the interface of $\pderiv{u}{x}$
in relation
to those of the local time of the process.

Finally, as one can easily check by matching the operators in (\ref
{1dLambdapde}) with the characterization of the infinitesimal operator
in (\ref{Admds}), the family of skew diffusion processes coincides with
the class of
Markov processes with scale and speed measures having piecewise
constant densities.

\begin{theorem}
Let $X$ be a regular diffusion process with speed measure $m$ and scale
measure $s$ having densities
%
%
\begin{eqnarray}\qquad
m'(x) &=& m^- \IND{(-\infty,0]}(x) + m^+
\IND{(0,\infty)}(x),\nonumber\\[-8pt]\\[-8pt]
s'(x) &=& s^- \IND{(-\infty,0]}(x) + s^+ \IND{(0,\infty)}(x)\nonumber
\end{eqnarray}
for some $m^+,m^-,s^+,s^- > 0$. Then $X$ is given by
%
%
\begin{eqnarray}
X = s_{\sqrt{D}}\bigl(B^{(\alpha)}\bigr)\qquad\qquad\qquad\nonumber\\[-8pt]\\[-8pt]
\eqntext{\displaystyle \mbox{with } D =
\frac{2}{m' s'}, \alpha= \frac{\sqrt{m^+ s^-}}{\sqrt{m^- s^+} + \sqrt{m^+
s^-}}.}
\end{eqnarray}
\end{theorem}

While there is no denying the importance of the conservative interface
condition in many physical applications, a primary goal of the present
article is to illustrate both the ubiquity and special effects of more
general interface conditions. This is especially relevant to certain
biological and ecological applications where awareness of such effects
might help to guide the determination of an appropriate interface condition.
For example, in the following section, we introduce notions of
\textit{natural occupation time} and \textit{natural local time}
for the processes arising in this more general class of models. These
are modifications
of the more standard mathematical definitions of occupation and local
time to adapt to
the physical units of the model, that is, so that occupation time is in
the units of time, for example.
An interesting consequence is that continuity and ordering properties
of these quantities can be obtained
that illustrate the effect of a particular interface condition in the
context of modeling
ecological and natural processes, for example, in relation to modeling
advection--dispersion of insect populations as considered in \citet
{cantrell2003spatial}, \citet{levin}.

\section{Applications in the Physical and Biological Sciences}\label{sec4}
\label{SectionApplications}

In this section several different areas of application are described.
It is in this section that
the Mathematics of Planet Earth theme is most clearly illustrated.
Each application area involves a distinct manifestation of an interface effect.

It is fitting to first note that the general interface conditions
introduced in
(\ref{1dLambdapde}) from a mathematical perspective
already arise naturally in a class of physical problems involving
heat conduction in heterogeneous media as follows. As treated, for
example, in
the classical reference \citet{CarslawJaeger},
the equation of conservation of thermal energy in a
thin rod composed of two semi-infinite rods with
heat capacity $\rho^{\pm}$ and heat conductivity $\kappa^{\pm}$,
respectively, is given by
%
%
\begin{equation}
\rho^{\pm} \,\frac{\partial u}{\partial t} = \kappa^{\pm}\,
\frac{\partial^2 u}{\partial x^2}
\end{equation}
with interface condition at $x=0$, given by
%
%
\begin{eqnarray}
u\bigl(t,0^+\bigr)&=&u\bigl(t,0^-\bigr),\nonumber\\[-8pt]\\[-8pt]
\kappa^+\,\frac{\partial u}{\partial x}\bigl(t,0^+
\bigr) &=& \kappa^-\,\frac{\partial u}{\partial
x}\bigl(t,0^-\bigr).\nonumber
\end{eqnarray}
In the notation used in this paper, $D^{\pm} = \frac{\kappa^{\pm}}{\rho
^{\pm}}$.
The fact that the interface condition only depends on the heat
conductivity coefficient
leads to the interface condition
%
%
\begin{equation}
\lambda\,\frac{\partial u}{\partial x}\bigl(t,0^+\bigr)
= (1-\lambda) \,\frac{\partial u}{\partial x}
\bigl(t,0^-\bigr),
\end{equation}
where $\lambda= \kappa^+/(\kappa^+ +\kappa^-)$.

The further collection of examples
provided below indicate various other contexts from biological,
environmental and physical sciences in
which the general interface conditions may arise.
We begin, however, with a return to an example of solute transport in
porous media for the first two illustrations of the theory. In
particular, the condition (\ref{conserveinterface}) applies at
the interfaces.

\subsection*{Heterogeneous Taylor--Aris Dispersion and Averaging Effects}
Taylor--Aris dispersion is well known throughout the physical and
biological sciences for its role in providing the effective rate of
spread of a solute immersed in a homogeneous fluid flow as given by
(\ref{adpde}) in the case of
Poisseuille flow directed along the horizontal axis of a cylindrical
tube ${\mathbf G} = [0,\infty)\times G$ in terms of the \textit{tube radius}
$R$ of the cross section $G$, the \textit{molecular diffusion
coefficient} $D$ and the \textit{maximum flow} $v_0$ (or \textit{cross-sectional average} $v_0/4$) of the parabolic flow
profile.\footnote{Motivated by considerations of the stability of a
viscous liquid to two-dimensional disturbances in a porous medium,
\citet{wooding} adapted their
analysis to obtain the corresponding formula for dispersion of a solute
in a unidirectional parabolic flow between two parallel planes
separated by a distance $R$. The geometry will effect the constants
appearing in the formulae for effective dispersion rates in ways that
are made clear by the general theorem of \citet{bhattgupta84}.}
In the case $v_0 = 0$ the
dispersion coincides with molecular diffusion, and when $D = 0$ the
dispersion of solute is aligned with the parabolic profile of the flow. The
relative contributions of these combined effects $(D > 0, v_0 > 0)$ are
captured time asymptotically in Taylor's remarkable insights, leading
to the celebrated formula for an \textit{effective dispersion rate}
%
%
\begin{equation}
\label{eqtaydisp} \Db= 2D + \frac{R^2v_0^2}{96D}.
\end{equation}
Although originally developed by \citet{taylor} and refined by
\citet{aris} using perturbation methods of partial differential
equations, this was subsequently shown to be a manifestation of the
central limit theorem for a concentration of Brownian motion
particles advected by the flow in \citet{bhattgupta84} for the
case of
Lipshitz continuous dispersion and drift coefficients. In this
context, the effective dispersion coefficient is a time-asymptotic
\textit{variance} parameter for the distribution of the position of an
immersed particle. In the presence of heterogeneity, as is currently
known, it had been loosely anticipated that the effective dispersion
would be modified by ``averaging;'' for example, see
\citet{gelharaxness1983}. In this section we will see that as a
result of an interface effect, the effective rate involves both \textit{arithmetic} and \textit{harmonic} averaging.

Consider (\ref{adpde}) in a cylindrical domain ${\mathbf G} = \R\times G$
with a cross-section $G \subset\R^d$, which is a closed interval in
the case $d=1$, or a bounded region with a smooth boundary if $d=2$.
Suppose the drift ${\mathbf v}$ is parallel to the $x_1$-axis and the
diffusivity is a diagonal matrix depending only on the transverse
variables. Namely, for $d=2$, ${\mathbf v} = (v_1,0,0)$ and $D =
\operatorname{diag}(D_1,D_2,D_3)$ with $v_1 = v_1(x_2,x_3)$ and $D_i =\break D_i(x_2,x_3)$
being positive and bounded away from zero, $i=1,2,3$. Let $c(t,\x)$ be
a solution, and consider its cross-sectional average,
%
%
\begin{equation}
\label{C} C(t,x):= \iint_{G} c(t,x,x_2,x_3)
\sud{x_2} \sud{x_3}.
\end{equation}
If $X(t) = (X_1(t),X_2(t),X_3(t))$, $t>0$ is the diffusion process
associated with the p.d.e. solved by $c$, then $C(t,\cdot)$ represents the
(nonnormalized) marginal distribution of the longitudinal coordinate
$X_1(t)$ for an initial uniform distribution
of the transverse coordinates $(X_2(t),X_3(t))$ on $G$.
The Taylor--Aris problem involves \textit{homogenized} parameters ${\bar
{v}}, \Db$ such that on large space--time scales $\lambda x$, $\lambda^2
t$, the weak limit
%
%
\begin{equation}
\label{longcenterscale}\quad \tilde{C}(t,x) \sud{x}:= \lim_{\lambda\rightarrow\infty}
C\bigl( \lambda^2 t, \lambda x + {\bar{v}} \lambda^2 t
\bigr) \lambda\sud{x}
\end{equation}
provides a centered solution
\[
\Cb(t,x) = \tilde{C}(t,x- {\bar{v}} t)
\]
to the homogenized partial differential equation,
%
%
\begin{equation}
\label{PDEC}\qquad \pderiv{\Cb} {t} = \frac{1}{2} \Db\,\pderiv[2]{\Cb} {x} - {
\bar{v}} \,\pderiv {\Cb} {x},\quad t \geq0, x \in\R.
\end{equation}

The homogenized parameters $\bar v$, $\bar D$ are in fact the result of
an ergodic theorem for the
transverse Markov process with reflecting boundary; see (\ref{bar}).
The following extension of \citet{bhattgupta84} can be obtained
for the case of a layered medium with piecewise continuous
coefficients; see \citet{ramirez2006} where the original idea
of \citet{bhattgupta84} to view the problem in accordance
with (\ref{longcenterscale}) as a functional central limit theorem for
$\{X_1(\lambda^2 t) - {\bar{v}} \lambda^2 t\}/\lambda$ as $\lambda\to
\infty$ is shown to carry over to piecewise continuous coefficients as well.

\begin{theorem}[(A generalized Taylor--Aris formula for piecewise
continuous coefficients)]
\label{GBF}
Assume\break  \mbox{$d=2$}. Let $\pi(\ud x_2 \sud x_3)$ be the uniform probability
measure on $G$,
and let $h$ be a solution in $L^2(G,\pi)$ to the boundary value problem
%
%
\begin{equation}\qquad
\cases{ \div(D_{2,3} \nabla h) = v_1 - {\bar{v}}, &
$(x_2,x_3) \in G$,
\vspace*{2pt}\cr
(D_{2,3} \nabla h)
\cdot{\mathbf n}_0 = 0, & $(x_2,x_3) \in\partial
G$,}
\end{equation}
where ${\mathbf n}_0$ denotes the outward normal vector of $G$ and $D_{2,3}
= \operatorname{diag}(D_2,D_3)$. Then, for any $t>0$, $x \in\R$,
and Borel measurable $A\subseteq\R$ with $|\partial A| = 0$,
%
%
\begin{eqnarray}
\label{limit}
&&\lim_{\lambda\rightarrow\infty} \int_A C\bigl(
\lambda^2 t,\lambda x + {\bar{v}} \lambda^2 t\bigr)
\lambda\sud{x} \nonumber\\[-8pt]\\[-8pt]
&&\quad= \int_A \Cb(t,x + {\bar{v}} t) \sud{x}\nonumber
\end{eqnarray}
with homogenized parameters
%
%
\begin{eqnarray}
\label{bar} {\bar{v}} &=& \iint_{G} v \pi(\ud
x_2 \sud x_3),\nonumber\\[-8pt]\\[-8pt]
\qquad\Db &=& \iint_{G}
\bigl\{D_1 + (D_{2,3} \nabla h) \cdot\nabla h\bigr\} \pi (\ud
x_2 \sud x_3).\nonumber
\end{eqnarray}
\end{theorem}

In the case $d=1$, $D_{2,3}$ is the scalar $D_2$ which is piecewise
constant, and we obtain
(see Figure~\ref{FigureTA})
the following corollary.

\begin{figure}[b]

\includegraphics{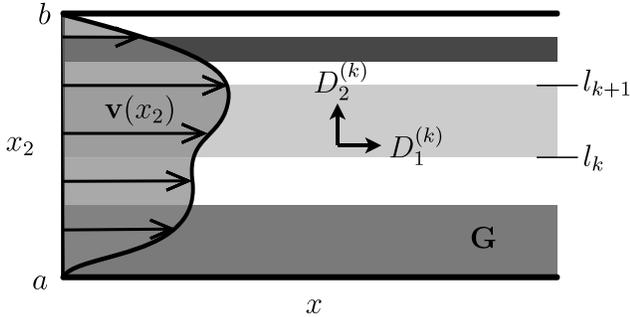}

\caption{Two-dimensional advection dispersion through a layered medium.
In the ongoing notation, $d=1$, $G = [a,b]$, ${\mathbf G} = \R\times
[a,b]$.} \label{layers}\label{FigureTA}
\end{figure}

\begin{cor}[(A generalized Taylor--Aris formula with piecewise constant
coefficients)]\label{ConstantLayers}
Assume $d=1$, $G=[a,b]$, and $D$ has the form
\begin{eqnarray*}
D &=& D(x_2) = \sum_{k=-m}^{M} D^{(k)}   \IND{[l_k,l_{k+1})}(x_2),\\
D^{(k)} &=& \lleft[
\matrix{
D_1^{(k)} & 0
\cr
0 & D_2^{(k)}}
\rright],
\end{eqnarray*}
where $a=l_{-m}<l_{-m+1}<\cdots<l_{M}<l_{M+1}=b$ is a collection of
interfaces partitioning $[a,b]$. If $D_1^{(k)}>0$ and $D_2^{(k)}>0$ for
all $k$,
then the limit (\ref{limit}) of Theorem~\ref{GBF} holds
with homogenized diffusion coefficient
%
%
\begin{eqnarray}
\label{bar2d} \Db&=& \sum_{k=-m}^{M} \biggl
\{D_1^{(k)} \frac{l_{k+1}-l_k}{b-a}\nonumber\\[-8pt]\\[-8pt]
&&\hspace*{28.6pt}{} + \frac{1}{D_2^{(k)}} \int
_{l_k}^{l_{k+1}} g(y)^2 \pi(\ud{y}) \biggr\},\nonumber
\end{eqnarray}
where $g$ is given by
%
%
\begin{equation}
\label{g} g(y) = \int_a^y
\bigl(v_1(x_2) - {\bar{v}} \bigr) \pi(
\ud{x_2}).
\end{equation}
\end{cor}
Thus, the first term of the effective dispersion rate (\ref{eqtaydisp}) is
replaced by a (weighted) arithmetic average, while the second term
involves a
(weighted) harmonic mean.
In particular, for $G = [-R,R]$ with a single interface
at $0$ separating media with diffusion coefficient $D^+$ and $D^-$,
respectively, and a parabolic velocity profile $v_1(x_2) = v_0
(1-(x_2/R)^{2})$, the formula is
%
%
\begin{equation}
D = D_a + \frac{4   v_0^2   R^2 }{ 945 D_h},
\end{equation}
where
%
%
\begin{equation}
D_a = \frac{D^+ + D^-}{ 2},\quad D_h = \frac{1 }{ {1}/{ D^+}
+ {1}/{D^-}}.\hspace*{-32pt}
\end{equation}

\subsection*{Physical Skew Diffusion and Stochastic Ordering of
Breakthrough Curves}
The topic addressed in this subsection
was originally initiated as a result of
observations resulting from
laboratory experiments designed to empirically test
and understand advection--dispersion in the presence of sharp
interfaces, for example,
experiments by \citet{kuo1999},
\citet{hoteit} and
\citet{Berkowitz09}.
Such laboratory experiments have
been rather sophisticated in the use of layers of sands and/or
glass beads of different granularities and modern measurement
technology. The specific interest is in
the effect of the interface condition on so-called breakthrough
curves, measuring the time required for
an injected concentration at one location to appear at another.
The basic phenomenon of interest to us here is captured by the following:

\textit{Question}. Suppose that a dilute solute is injected at a
point $y$ units to the left of
an interface at the origin and retrieved at a point $y$ units to the
right of the interface. Let $D^-$ denote
the (constant) dispersion coefficient to the left of the origin and
$D^+$ that to the right, with,
say, $D^- < D^+$ (see Figure~\ref{FigureInterfacial}). Conversely,
suppose the solute is injected at a point $y$ units to the right of the
interface and retrieved at a point $y$
units to the left. \textit{In which of these two symmetric arrangements
will the immersed solute most rapidly break through at the opposite
end}?

\begin{figure}

\includegraphics{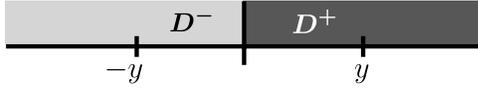}

\caption{Interfacial schematic.}\label{FigureInterfacial}
\end{figure}

The following results indicate that the question above can be answered
by investigating the asymmetries in the hitting times of skew Brownian
motion and skew diffusion.
\begin{lemma}
\label{skewtoelastic}
Fix $y \ge0$ and let $H^{(\alpha)}_y = \inf\{t\ge0\dvtx  B^{(\alpha)}(t) =
y\}$.
If $0< \alpha<1/2$, then
\begin{eqnarray*}
\P_{-y}\bigl(H_{y}^{(\alpha)} > t\bigr) &<& \frac{\alpha}{ 1-\alpha} \P
_y\bigl(H_{-y}^{(\alpha)} > t\bigr)\\
&<& \P_y\bigl(H_{-y}^{(\alpha)} > t\bigr),\quad
t > 0.
\end{eqnarray*}
\end{lemma}

Recall that rescaling space by the respective diffusivities symmetrizes
the transition probabilities
when $\alpha= \alpha^*$, that is, for physical skew diffusion $X^*(t)
= s_{\sqrt{D}}(B^{(\alpha^*)}(t))$.
The following stochastic ordering of first passage times for physical
skew diffusion provides a simple probabilistic basis for
the symmetries and asymmetries
predicted in experiments cited\break  above. The proof is by a coupling and
relies on an interesting balance between
the specific conservative transmission parameter $\alpha^*$, as well
as the respective scalings on either side of the interface; see
Appuhamillage et~al. (\citeyear{AppuhamAAPCorrect,AppuhamilageAAP}) for details.
%
\begin{theorem}
\label{thmasympassagetime}
Let $H^*_y = \inf\{t\ge0\dvtx  X^*(t) = y\}$. Then,
for $y > 0$ and $D^- < D^+$,
\begin{eqnarray*}
\P_{-y}\bigl(H^*_y > t\bigr) &\leq&\frac{\sqrt{D^-}}{\sqrt{D^+}}
\P_{y}\bigl(H^*_{-y} > t\bigr) \\
&<& \P_{y}
\bigl(H^*_{-y} > t\bigr),\quad t\ge0.
\end{eqnarray*}
\end{theorem}

To gain an alternative perspective on this phenomena, one may
compute and compare the concentration curves as a function $t\to
u(t,y)$ and $t\to u(t,-y)$
for a point injection at the interface; see
Appuhamillage et~al. (\citeyear{AppuhamAAPCorrect,AppuhamilageAAP}).

\subsection*{Interfacial Effects on Natural Residence and Local Times
of Skew Diffusions}

Within the ecology literature there is a recognition of the role of
interfaces in ``directing'' movement from one habitat to another (e.g.,
see Fagan, Cantrell and Cosner, \citeyear{fagan1}; \citecs{cantrell2003spatial}, page 112, and numerous
references therein, as well as \citecs{levin}, page 265). The main
point of the example described here is to highlight a natural role of
nonconservative interfacial conditions in the models involving insect
dispersal. Specifically, we will examine the effect of the interface on
functionals such as residence times.

Fender's Blue butterfly provides a specific example that has been
analyzed fairly extensively in both field experiments as well as mathematics.
Fender's Blue is an endangered species of butterfly found
in the Pacific Northwestern United States. The primary habitat patch is
Kinkaid's Lupin flower.

Ecologists have focused substantial fieldwork efforts in examining the
way in which organisms respond to habitat edges and the relationship to
residence times in Lupin patches; see \citet{Schultz2001}.
Sufficiently long residence (occupation) times
are required for pollination, eggs, larvae and ultimate sustainability
of the population. Empirical evidence points to a skewness in random
walk models for butterfly movement at the path boundaries that have led
to ``biased random
walk'' and skew Brownian motion models in
\citet{Schultz2001}, \citet{cantrell2003spatial}, \citet
{levin}, \citet{fagan1}, \citet{ovas}. The
determination of proper interface conditions is primarily a statistical
problem in such applications. However, as illustrated below, the role
of local interfacial conditions is reflected in the behavior of
residence times
in ways that may be useful to the identification of interface conditions.
In the framework of one-dimensional advection--dispersion one is
therefore lead to consider the interface conditions
(\ref{1dLambdapde}) generalizing the conservative interface
condition
(\ref{conserveinterface}).

Note that $\lambda=0,
\lambda= 1$ correspond to reflection at the interface,
while $\lambda= \frac{D^+}{D^+ + D^-}$ is the
conservative interface condition (\ref{conserveinterface}) that gives
rise to the process $X^*$.
In particular, at the extremes the residence times of the
positive half-line are obviously quite distinct.
The following result interpolates between these extremes.
The proof
exploits the basic property of skew
Brownian motion noted at the outset in (\ref{skeweffect}),
and essentially that
%
%
\begin{equation}
\label{posprob} \P_0\bigl(B^{(\alpha)}(t) > 0\bigr) = \alpha,\quad
t > 0.
\end{equation}
This is easily checked from the
definition and, intuitively, reflects the property that the
excursion interval $J_{n(t)}$ of $|B|$
containing $t$ results in a $[A_{n(t)} = +1]$ coin flip with probability
$\alpha$.

The following theorem involves a modification of the usual mathematical
definition of \textit{occupation time}, for example, as given in standard
references such as \citet{revuzyor}, in that integration is with respect
to the Lebesgue measure in place of quadratic variation. We refer to
this modification as
\textit{natural occupation time}.

\begin{definition}
\label{natoccuptime}
Let $X$ be a continuous semimartingale.
The \textit{natural occupation time} of a Borel set $A$ by
$X$ in time $[0,t]$ is defined by
\[
\tilde{\Gamma}^X(A,t) = \int_0^t
\IND{A}\bigl(X(s)\bigr) \sud s.
\]
\end{definition}

One
may note that this modification puts occupation time in the natural
units of ``time,'' while mathematical local time is in units of (area)
``spatial length squared.'' As such,
natural occupation time seems to be the more appropriate representation
of \textit{residence time}
measurements, and we use it here for identifying regularities and
properties of interest to the
applications.
Mathematical occupation time, on the other hand, has important roles to
play in other theoretical contexts.

\begin{theorem}
\label{thmoccuptime}
Let $X^{(\alpha(\lambda))}$ denote skew diffusion defined in (\ref
{defXLambda}) for the dispersion
coefficients $D^+, D^-$ and interface parameter $\lambda$. Denote
natural occupation time
processes by
\[
\tilde{\Gamma}_\lambda^+(t) = \int_0^t
\IND{(0,\infty)}\bigl(X^{(\alpha
(\lambda))}(s)\bigr) \sud s,\quad t\ge0.
\]
Similarly, let $\tilde{\Gamma}_\lambda^-(t) = t- \tilde{\Gamma}_\lambda
^+(t), t\ge0$. Then,
\[
\EXP\bigl(\tilde{\Gamma}_\lambda^+(t)\bigr) > \EXP\bigl(\tilde{
\Gamma}_\lambda^-(t)\bigr),\quad t>0,
\]
if and only if
\[
\lambda>
{\sqrt{D^+}\over\sqrt{D^+}
+\sqrt{D^-}}
\]
with equality when $\lambda= {\sqrt{D^+}\over\sqrt{D^+} +\sqrt{D^-}}$.
\end{theorem}

It is noteworthy, therefore, that under the conservative interface condition
more time is spent in the more volatile habitat, making such models
questionable for many ecological contexts involving animal dispersion.\vadjust{\goodbreak}

The conservative interface condition can also be characterized as the
unique interface
condition that gives the continuity of an analogous \textit{natural local
time} defined as follows.

\begin{definition}
\label{natloctime}
Let $X$ be a continuous semimartingale. The \textit{natural local time at $a$}
$\tilde{L}^X(t,a) = \frac{1}{2}(\tilde{L}^{X,+}(t,a) + \tilde{L}^{X,-}(t,a))$
of $X$ is defined by
\[
\tilde{L}^{X,+}(t,a) = \lim_{\varepsilon\downarrow0}
{1\over\varepsilon}\int_0^t\IND{[a, a+
\varepsilon)}\bigl(X(s)\bigr)\sud s
\]
and
\[
\tilde{L}^{X,-}(t,a) = \lim_{\varepsilon\downarrow0}
{1\over\varepsilon}\int_0^t\IND{(a-
\varepsilon, a)}\bigl(X(s)\bigr) \sud s,
\]
provided that the indicated limits exist almost surely.
\end{definition}

With this definition, one has that
%
%
\begin{equation}
\tilde{\Gamma}^X(A,t) = \int_A
\tilde{L}^X(t,\ud a)
\end{equation}
in complete analogy with the standard relation between local time and
occupation time defined
using the quadratic variation of the process $X$.

Recall that in the particular case of skew Brownian motion, the
quadratic variation is simply\break  $\langle B^{(\alpha)} \rangle(t) = t$. Therefore,
the ``symmetric local time'' $\frac{1}{2}(L(^X(t,a) + L^X(t,a^-))$ [see
\citecs{revuzyor} and equation (\ref{altDefLeft})] agrees with the
natural local time just defined. Moreover, the following relations
among one-sided and symmetric local times at $0$ are known; for
example, see \citet{Ouknine}:
%
%
\begin{eqnarray}
\label{localtimerelations} 2\alpha\tilde{L}^{B^{(\alpha)},+}(t,0) &=&
\tilde{L}^{B^{(\alpha)}}(t,0),\nonumber\\[-8pt]\\[-8pt]
2(1-\alpha) \tilde{L}^{B^{(\alpha)},-}(t,0) &=&
\tilde{L}^{B^{(\alpha)}}(t,0).\nonumber
\end{eqnarray}
In particular, 
the symmetric (natural) local time is continuous if and only if
$\alpha=1/2$.

The next theorem, a version of which was originally developed in
\citet{Appuhamillage2012}, extends the continuity of natural local
time to the more general framework of the present paper. While the
purpose here is not to explore the generality for which natural local
time exists among all continuous semimartingales, according to the
following theorem
it does exist for skew diffusion. Moreover, continuity has a special
significance for
the determination of parameters.

\begin{figure*}

\includegraphics{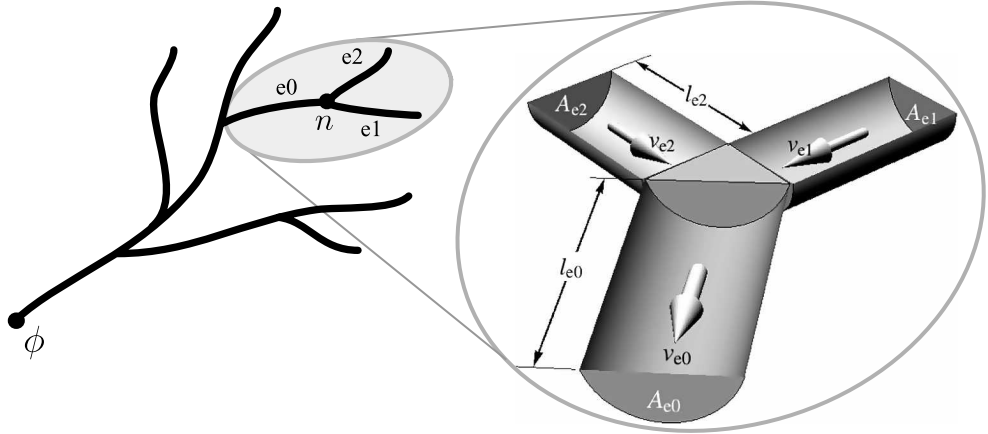}

\caption{Schematic of a river network $\Gamma$ with root node $\phi$.
The inset shows three edges connected at an internal node $n$
representing the junction where two tributaries merge to form a new
channel.}\label{FigureNetworkY}
\end{figure*}

\begin{theorem}
\label{natloctimeofnatdiff}
Let $X^{(\alpha(\lambda))}$ be the skew diffusion process with
parameters $D^{\pm}, \lambda$.
Then the natural
local time of $X^{(\alpha(\lambda))}$ at $0$
is continuous if and only if
$\lambda= {D^+\over D^+ +D^-}$, that is, if and only if $\alpha
(\lambda)=\alpha^*$
and, thus,
$X^{(\alpha^*)}$ is the physical skew diffusion.
\end{theorem}

Thus, while at the macroscale of deterministic particle concentrations
the determination of the transmission parameter $\alpha^*$ may
be viewed as a consequence of the continuity of flux at the interface,
at the scale of stochastic
particle motions, it is determined by a condition of continuity of
natural local time at the interface.

\subsection*{Dispersion of Organisms in River Networks}

River networks are known to control the flux of water and sediment over
most landscapes on the planet earth. Moreover, transport of water,
organisms, sediment, \mbox{nutrients} and contaminants on river networks plays
a central role in modern hydrology and ecology. River networks
constitute, in particular, fundamental ecosystems whose populations are
dependent upon the interconnectivity and heterogeneity of the different
reaches that form the network
\citep{fagan2}.

Mathematically, river networks are modeled as directed binary graphs. A
long tradition of research in hydrology and geomorphology has narrowed
the class of graphs observed in natural river basins, the
relationships between physical variables involved in transport and the
topological properties of such networks; for example, see
\citet{RodriguezIturbe2001}, \citet{peckham},
\citeauthor{barndorff}\break
(\citeyear{barndorff}) and
references therein. It is therefore natural to extend the linear
advection--diffusion (\ref{1dskewpde}) to a binary graph in an effort to
advance the understanding of the relationship between network topology,
physical properties of rivers and dispersal of organisms.

The first steps toward an extension of skew Brownian motion to an
infinite star-shaped graph was introduced by \citet{walsh1} as a natural
mathematical extension of skew Brownian motion on $\R$. A general
theory of advection--diffusion processes on arbitrary graphs was
subsequently
initiated by
\citet{Freidlin1993}.

To fix ideas in the context of river networks, consider a connected
binary directed tree graph $\Gamma$ as depicted in Figure~\ref{FigureNetworkY}. Each edge $e$ models a stream reach of length $l_e$
between two junctions and is assumed to be isomorphic to the interval
$[0,l_e]$. Also, each edge $e$ has associated strictly positive
parameters $v_e$, $A_e$ and~$D_e$, denoting the mean water velocity,
cross-sectional area and diffusion coefficient of the organisms in that
reach. The endpoints $x=0$ and $x=l_e$ correspond to the downstream and
upstream nodes, respectively. The set of nodes in $\Gamma$ can be
divided into three subsets: the singleton root node $\phi$, the set
$I(\Gamma)$ of internal nodes connecting three edges, and the set
$U(\Gamma)$ of upstream nodes $n$ of ``tributary edges'' or ``leaves'' of
$\Gamma$.

Considering the spatio-temporal evolution of the density of suspended
organisms in $\Gamma$, and imposing conservation of mass throughout,
one arrives at the following extension of (\ref{1dskewpde})
to the network:
%
%
\begin{eqnarray}
\label{treeskewpde} \pderiv{u_e} {t} = \frac{1}{2}
D_e \,\pderiv[2]{u_e} {x} - v_e \,\pderiv
{u_e} {x},\nonumber\\[-8pt]\\[-8pt]
\eqntext{x \in[0,l_e], e \in\Gamma,}
\end{eqnarray}
where $u_e$ denotes the restriction of the function $u$ to edge $e$.
Only functions that are continuous on each edge of $\Gamma$ and twice
continuously differentiable on the interior of each edge are
considered. For an internal node where edges $e1$, $e2$ join to form
edge $e0$, the appropriate extension of (\ref{conserveinterface}) reads
%
%
\begin{eqnarray}
\label{treeconservativeinterface}\qquad
&& A_{e0} D_{e0}\,
\pderiv{u_{e0}} {x}(t,l_{e0}) \nonumber\\[-8pt]\\[-8pt]
&&\quad= A_{e1}
D_{e1} \,\pderiv{u_{e1}} {x}(t,0) + A_{e2}
D_{e2} \,\pderiv{u_{e2}} {x}(t,0).\nonumber
\end{eqnarray}
Here, we also have assumed that water discharge is conserved at river
junctions, namely, $A_{e1} v_{e1} + A_{e2} v_{e2} = A_{e0} v_{e0}$.

Several different behaviors can be prescribed at the boundary nodes of
$\Gamma$. In particular, following \citet{Speirs2001kx},
\citet{lutscheretal2005}, \citet{lutscheretal2006}, one may
consider an ecological
scenario where organisms do not leave the network through channel
sources, and an abrupt change of flow conditions occur at $\phi$ that
removes organisms from $\Gamma$, for example, a waterfall, a fast
flowing river, a lake, the ocean or human disturbances. This can be
coded mathematically by requiring
%
%
\begin{eqnarray}
\label{treebdrynodes} u(t, \phi) &=&0,\nonumber\\[-8pt]\\[-8pt]
\pderiv{u} {x}(t,n) &=& 0,\quad n \in U(\Gamma), t
\geq0.\nonumber
\end{eqnarray}

As shown in \citet{Freidlin00}, \citet{Freidlin1993}, the
spatial operator on the left-hand side of (\ref{treeskewpde}) along
with conditions (\ref{treeconservativeinterface}, \ref{treebdrynodes})
is the infinitesimal generator of a Feller Markov process $X = \{X(t),
t \geq0\}$ on $\Gamma$ with continuous sample paths that can be
written as $X(t) = (x(t),e(t))$ with $e(t)$ being the edge the process
occupies at time $t$, and $x(t) \in [0,l_{e(t)}]$, $t \geq0$.
Moreover, one has the analogous representation to (\ref{skewDiffSDE}):
there exists a one-dimensional Brownian motion $B$ and an increasing
process $L$ such that
%
%
\begin{equation}\qquad
\ud x(t) = \sqrt{D_{e(t)}} \sud B(t) - v_{e(t)} \sud t + \ud
L(t),
\end{equation}
where $L$ only increases when $x(t) = 0$. The three-way heterogeneity
at internal nodes has a skewing effect on the sample paths analogous to
property (\ref{skewnessX}) of skew diffusion: let $H^x_\eps= \inf\{t
\geq0\dvtx  x(t) = \eps\}$ and $n$ denote the node connecting edges
$e0,e1,e2$, then
%
%
\begin{eqnarray}&&
\lim_{\eps\to0^+} \P\bigl(e\bigl(H_\eps^x
\bigr) = ei | X(0) = n\bigr) \nonumber\\
&&\quad= \frac
{A_{ei}D_{ei}}{A_{e0}D_{e0}+A_{e1} D_{e1}+A_{e2}D_{e2}},\\
\eqntext{i =
0,1,2.}
\end{eqnarray}

An important contribution of advection--diffusion models in riverine
ecology revolves about the classical ``drift paradox,'' whereby it was
observed in \citet{muller1954investigations} that although individual
organisms in streams are subject to downstream drift, the average
location of the population is not observed to move downstream over
time, and thereby persists. In this regard,
\citet{Speirs2001kx}, \citet{lutscheretal2005}, \citet
{lutscheretal2006} obtain
conditions on the channel length, drift velocity and population
dynamics under which the population as a whole can persist along a
single channel assuming that the movement of individuals is given by an
advection--diffusion process of the form (\ref{adsde}). Results of this
type define resolutions of the drift paradox that can be useful to
managers of ecological preserves.

Ecological persistence problems on river networks involve models in
which individuals move within $\Gamma$ via a jump process: an organism
initially located at $y \in\Gamma$ jumps to the position $X(\tau_\sigma
)$, where $\tau_\sigma$ is an exponentially distributed random variable
that represents the time the individual spends dispersing within the
water column. The resulting \textit{dispersal kernel}, as it is known in
the ecological literature, is therefore given by
%
%
\begin{eqnarray}
\label{dispersalkernel} k(y,x) &=& \P_y\bigl(X(\tau_\sigma)
\in\ud x\bigr) \nonumber\\
&=& \int_0^\infty\sigma
e^{-\sigma t} p(t,y,x) \sud t,\\
\eqntext{x,y \in\Gamma,}
\end{eqnarray}
where $p(t,y,x)$ are the transition probability densities of $X$, the
fundamental solution to problem (\ref{treeskewpde}), (\ref
{treeconservativeinterface}) and (\ref{treebdrynodes}).

In the solutions noted above by \citet{Speirs2001kx},
\citet{lutscheretal2005}, \citet{lutscheretal2006}, the
evolution of the population $u$ in the single channel $[0,l]$ is assumed
to be given by\footnote{Although \citet{lutscheretal2005} view
(\ref{lutschereqn}) as
``derived'' from a regime-switching model, the argument is flawed. On
the other hand, one may simply view (\ref{lutschereqn}) as a distinct
population model.
\citet{felderway} have recently shown, however, that the
conditions on the parameters for persistence for the regime-switching
model differ from those for
(\ref{lutschereqn}).}
%
%
\begin{eqnarray}
\label{lutschereqn} {\partial u\over \partial t} (t,x) &=& ru(t,x) - \lambda
u(t,x)\nonumber\\[-8pt]\\[-8pt]
&&{} +
\int_0^l \lambda k(y,x)u(t,y) \sud y,\nonumber
\end{eqnarray}
where $r > 0$ is the net population growth rate\footnote{More general
dynamics are considered, however, the results depend on the linear form
in (\ref{lutschereqn}).} and $\lambda$ is the jump rate. Of
course, in the case of a single channel there is not an interface and
the model for
$p(t,y,x)$ is simply Brownian motion with drift $v$ and diffusion
coefficient $D$.\footnote{\citet{Speirs2001kx},
\citet{lutscheretal2006} permit more general models for
$p(t,y,x)$, although the Brownian motion model is a primary example.}
Persistence is defined by instability of the solution $u = 0$ to (\ref
{lutschereqn}).

An extension from the interval $[0,l]$ to tree networks $\Gamma$ was
developed in \citet{ramirez2012green}, wherein the dispersal kernel
(\ref{dispersalkernel}) is explicitly solved. Moreover,
\citet{ramirez2012green} permits nontrivial events of upstream
migration, which are proposed in \citet{Speirs2001kx} to be the
key to
explaining the drift paradox. In particular,
\citet{ramirez2011population} provides bounds on the minimum net growth
rate $r$ of individuals required for persistence of the population
evolution in $\Gamma$ via the dispersal kernel
(\ref{dispersalkernel}).

\subsection*{Coastal Upwelling, Fisheries and Continuity of Natural
Local Time}

Night satellite images of the earth show a striking concentration of
fishing flotillas exploiting the ocean bounty
off the coast of southern South America between approximately latitude
40 and 50 degrees south. This activity takes place in a narrow strip
that follows the continental shelf break of South America where the cold
nutrient rich waters of the Malvinas current reach the surface of the
Atlantic Ocean in a process described as upwelling (see Figure~\ref{FigureMalvinas}).
A mathematical
%
\begin{figure}

\includegraphics{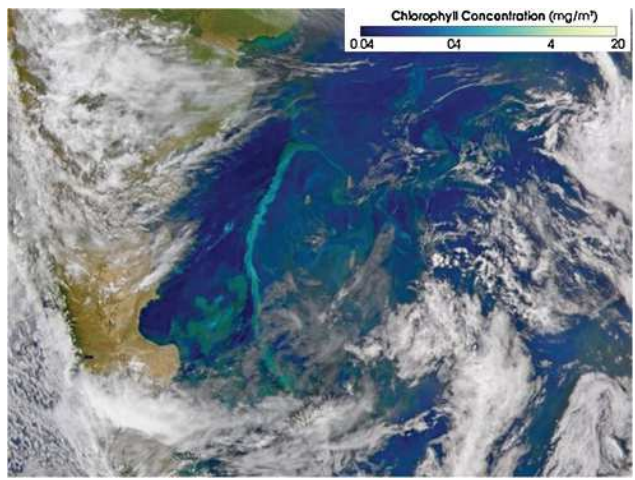}

\caption{Phytoplankton bloom in Malvinas/Falklands current off the
Atlantic coast of southern South America. Provided by the SeaWiFS
Project, NASA/Goddard Space Flight Center and
\mbox{ORBIMAGE}.}\label{FigureMalvinas}
\end{figure}
model that approximately describes the location of the upwelling is
given by the \textit{arrested topographic wave} equation. Obtained under
various simplifying assumptions, such as hydrostatic approximation and
geostrophic balance, the equations determine the ocean-free surface
elevation $\eta(x,y)$ as the solution of
%
%
\begin{equation}
\label{ATW} \frac{\partial\eta}{\partial y} =
-\frac{r}{f}\frac{1}{h(x)}\,
\frac
{\partial^2\eta}{\partial x^2}.
\end{equation}
Here $x$ is the distance from the shore, $y$ is the along-shore coordinate,
$r>0$ is the bottom friction, $f<0$ is the Coriolis parameter and
$h(x)$ is the
derivative of the depth of the ocean at $x$. As part of the derivation
of the
equation, the orientation of the along-shore axis $y$ is determined by
the direction of
motion of the current; see \citet{MatanoPalma} and references within.
As a consequence, (\ref{ATW}) can be thought of as a diffusion equation
where $y$ plays the role of time.

The main features of the upwelling process are obtained by considering
a bottom
topography in which the slope $h(x)$ is piecewise constant, $h^{\pm}$
with a
discontinuity at the continental shelf brake $x=l$. At the
``interface''
the conservation of mass transport by the current leads to the conditions
%
%
\begin{eqnarray}
\label{interfaceATW} \eta\bigl(l^-,y\bigr) &=&
\eta\bigl(l^+,y\bigr),\nonumber\\[-8pt]\\[-8pt]
\frac{\partial\eta}{\partial x}\bigl(l^-,y\bigr) &=&
\frac{\partial\eta}{\partial x}\bigl(l^+,y\bigr).\nonumber
\end{eqnarray}
This equation corresponds to (\ref{1dLambdapde}) with
\[
D^\pm= -\frac{r}{f h^\pm},\quad \lambda=1/2.
\]
While the conservative interface conditions (\ref{conserveinterface})
are ubiquitous in the hydrological literature, the conditions
(\ref{interfaceATW}) are mathematically natural. If one follows the
general theory of time changes in the context of martingale problems
as presented in \citet{StrookVaradhan}, the interface conditions of
the arrested topographic wave are the ones that can obtained by a
direct application of this theory. Indeed, with the time change
$\tau(t) = \int_0^t \frac{1}{D(w(s))} \sud s$
one obtains that $X(t)= B({\tau(t)}), t\ge0$, is the Markov process
corresponding to the problem (\ref{ATW})
with interface condition (\ref{interfaceATW}).
Alternatively, one can obtain $X(t) = s_{\sqrt{D}}(B^{\alpha^{\#}}(t)),
t\ge0$,
where $\alpha^{\#} = \sqrt{D^-}/\break (\sqrt{D^+} + \sqrt{D^-})$.

\section{Complementary Results and Open Problems}\label{sec5}

The main goal of this article has been to develop basic pathways to the
frontiers of advection--disper-\break sion research in the presence of
interfacial effects from a probabilistic point of view. The example of
one-dimensional processes with point interfaces is rich enough to
provide realistic illustrations of the diverse effects on quantities
arising in the applied sciences and engineering, however, it falls
far short of a general mathematical framework. In addition, even in the
one-dimensional context,
the examples were selected to highlight various significant interfacial
effects, but without an attempt to be comprehensive. However, the
relative consequences of these effects
do not seem to be widely recognized in the science literature in terms
of the specification of the interfacial condition. As a result, the
conservative interface condition is often adopted as the default consideration.

In this section we indicate some related results and open problems at
the frontiers of research in this general area of advection--dispersion.

An alternative
approach has also been partially developed by \citet{portenko1},
also see
\citet{portenko}, in the context of pdes whose coefficients may be
generalized
functions. Specifically, \citet{portenko} permits singular drift
terms but requires smooth
dispersion coefficients. A~somewhat heuristic development of ideas
along these lines in the
context of dispersion in porous media was explored in \citet
{labolle} and \citet{gravner} that may also provide effective
approaches to problems of this type. Certainly this provides an
intriguing mathematical framework to explore, especially for problems
in\break  higher dimensions.

The definitions pertaining to breakthrough curves have various not
necessarily equivalent formulations in the science and engineering
literature. While first passage time is of obvious probabilistic
interest, in the presence of advection the time profile of the
\textit{flux-averaged concentration} at a point is also sometimes
adopted. The flux-averaged concentration is expressed in terms of the
spatial derivatives of $u(t,x)$; see \citet{AppuhamilageWRR}. For
the case of a discontinuous medium, this means that the concentration
at the interface depends upon the derivatives on both the left and the
right at the interface. However, these gradients evolve differently in
the coarse and fine media. This means that, at the interface, the
gradients of the concentration do depend upon the configuration
(fine-to-coarse versus coarse-to-fine), and this dependency essentially
breaks the symmetry that can be observed for the breakthrough curve of
the resident concentration. This provides an alternative response to
the fine-to-coarse vs. coarse-to-fine breakthrough curves, as this is
explicitly computable in the presence of advection (see
Appuhamillage et~al., \citeyear{AppuhamilageWRR,AppuhamilageAAP}).
The determination of
the corresponding first passage times is unsolved in this generality.
However, \citeauthor{appsheldon}\break  (\citeyear{appsheldon}) recently computed an explicit
expression for the density of the first passage time of skew Brownian
motion.

The numerical methods by \citet{martineztalay} and \citet
{bokil2} are
more generally valid for piecewise continuous in place of piecewise
constant diffusion coefficients. However, the methods are exclusively
for the one-dimensional problems. The corresponding problems in a
higher dimension are largely untreated. Similarly, as suggested by the
examples, in many applications of biological interest in which the
conservative interface condition is inappropriate, the determination of
the proper interface condition involves statistical inference. The
papers by \citet{brillinger} and \citet{brill} provide some
perspective
on statistical inference for stochastic differential equations with
smooth coefficients with ecological applications that might be expanded
to this context. The upwelling example treated here is also a resource
for possible explanation of migratory patterns being monitored by ocean
ecologists in the context of oceanic flow properties; for example, see
\citet{acha2004marine}.

The three theorems presented in Section~\ref{sec3} illustrate different
approaches to deal with the effects of discontinuities of
coefficients in the model equations. As problems in earth,
physical and social sciences present situations in which abrupt
changes in the parameters are more common, there is a need
to develop an understanding from diverse points of view of the
effect of these discontinuities in easily observable (macroscopic)
quantities, as well as in refined, local properties of the processes
involved. The results in Section~\ref{sec4} illustrate some of these effects
in local time and occupation time as well as in the sustainability of
ecosystems.

As remarked in previous sections, the tools that have been developed so
far apply to one-dimensional problems, or problems on networks that
preserve a one-dimensional structure along its branches. An
outstanding open problem is to develop a more comprehensive approach to
problems with discontinuities across hypersurfaces in several space
dimensions. In this direction, \citet{Peskir2007} obtains an extension
of the It\^{o}--Tanaka formula to the case where discontinuities occur
across the graph of a function. This result, however, does not apply to
the case when discontinuities occur, for example, across the surface of
a sphere and severely limits its applicability since in many problems
of interest, the discontinuities occur across the boundary of bounded
sets. Two cases have been investigated in this regard:
\citet{decamps2006} constructs a skew Bessel process from its
scale and
speed measures and proves that such a process is the radial component
of a ``generalized diffusion process'' (in the language of
\citet{portenko1}) with a drift term concentrated on the boundary
of a
sphere. Second, the family of planar diffusions with rank-dependent
diffusion coefficients, thoroughly studied in
\citet{fernholz2011planar}, \citet{Fernholz2012fk}, include
the case of a
diffusion process in $\R^2$ with discontinuous diffusion coefficient
along the line $x=y$.

While full generalization of skew diffusion to problems in dimension
greater than one is yet to be completed, some progress has been made
for heterogeneous diffusion on graphs, as seen, for example, in Section~\ref{SectionApplications} in the context of river networks.
Additionally, for strong-swimming species whose movements are not
dependent on the water velocity, \citet{Gutierrez2011be} have used
analysis on the operator in (\ref{treeskewpde}) with $v_e =0$ for all
edges to study the time required to eradicate invasive species in a
river network.

At a more foundational level, \citet{hairer1} obtain a one-dimensional
skew diffusion as a limit in the diffusive scale of solutions of
stochastic differential equations with a periodic drift coefficient
outside an interval centered at the origin. In the limit the diffusion
coefficients are determined in a classical manner, while the skewness
is characterized in terms of a Zvonkin type transform of the drift.
Moreover, in \citet{hairer2} a similar limit is analyzed for a
stochastic differential equation with periodic drift in all directions
outside a finite region centered on a hyperplane in $\R^k$. The
limiting diffusion has an infinitesimal generator with discontinuous
coefficients for which the diffusion coefficient is classically
determined. In turn, the domain of the generator determining the
interface conditions is characterized through relations on the normal
derivatives from both sides of the hyperplane and tangential
derivatives. Such results illustrate the promise of a rich
mathematical theory as research on interfacial effects goes forward.

A perhaps intermediate step between skew diffusion and diffusion on
graphs is the case of problem (\ref{1dskewpde}) with a piecewise
diffusion coefficient taking more than two or an infinite number of
values. Namely, one might consider, as in Corollary \ref
{ConstantLayers}, a~one-dimensional medium with multiple interfaces. Define
%
%
\begin{equation}
\label{Dmultiple} D(x) = \sum_{k \in\Z} D^{(k)}
\IND{(x_k,x_{k+1})}(x),
\end{equation}
where $\{x_k\dvtx  k \in\Z\}$ is a sequence of real numbers with no
accumulation points and $\{D^{(k)}\dvtx  k \in\Z\}$ are positive numbers
uniformly bounded away from zero. The flux continuity condition (\ref
{conserveinterface}) extends to the multiple interface case as
$D^{(k-1)} f'(x_k^-) = D^{(k)} f'(x_k^+)$, $k \in\Z$. The associated
process is denoted by $X_M^*$ and is referred to as ``multiple skew
diffusion.''

Using a framework very similar to that presented in Section~\ref{SectionAlternatives}, \citet{ramirez2011multi} characterizes
$X^*_M$ in terms of a process exhibiting skewness of paths at several
points: ``multiple skew Brownian motion.'' More precisely, let $\ab=
\{\alpha_k\dvtx  k \in\Z\}$ be a sequence with $\alpha_k \in
(0,\frac{1}{2}) \cup(\frac{1}{2},1)$ and consider interfaces $\{y_k\dvtx\break   k
\in\Z\}$ with no accumulation points. The process $B^{\ab}$ is defined
as the Feller process with generator $\A^{\ab} f = \frac{1}{2} f''$
acting on functions $f \in C_b(\R)$ that are twice continuously
differentiable inside each $(y_k,y_{k+1})$ and that satisfy $\alpha_k
f'(y_k^+) = (1-\alpha_k) f'(y_k^-)$ for all $k \in\Z$. Not
surprisingly, ``multiple skew diffusion'' is given by
%
%
\begin{equation}
X_M^*(t) = s_{\sqrt{D}} \bigl( B^{\ab^*} (t)\bigr),
\end{equation}
where $s_{\sqrt{D}}$ denotes here the continuous piecewise linear
function with $s_{\sqrt{D}}(0) = 0$ and $s_{\sqrt{D}}'(x) = D^{(k)}$
for $x \in(x_k,x_{k+1})$, and $B^{\ab^*}$ is multiple skew Brownian
motion with interfaces $\{s_{\sqrt{D}}^{-1}(x_k)\dvtx  k \in\Z\}$ and
skewness sequence
%
%
\begin{equation}
\alpha_k = \frac{\sqrt{D_k}}{\sqrt{D_k} + \sqrt{D_{k-1}}},\quad k \in \Z.
\end{equation}

Processes of the form $X_M^*$ can be used in the context of transport
within layered media, as exemplified in Figure~\ref{FigureTA}. In
particular, \citet{ramirez2011multi} analyzes such a process to
determine the asymptotic behavior of particles undergoing
advection--diffusion on a periodic infinite layered medium composed of
two phases: a matrix of slow diffusive transport with periodic cracks
where fast diffusion dominates.

\section*{Acknowledgments}

The authors are especially grateful for the careful reading and comments
by the reviewers that greatly improved the exposition. Professor
Enrique Thomann acknowledges fellowship support of the NSF Institute
for Mathematics and its Applications
at the University of Minnesota during the academic year 2012--2013.
The authors Thomann and Waymire were supported in part by NSF Grants
DMS-11-22699 and DMS-10-31251.




\begin{thebibliography}{73}

\bibitem[\protect\citeauthoryear{Acha et~al.}{2004}]{acha2004marine}
\begin{barticle}[auto:STB|2013/09/19|12:14:10]
\bauthor{\bsnm{Acha},~\bfnm{E.~M.}\binits{E.~M.}},
  \bauthor{\bsnm{Mianzan},~\bfnm{H.~W.}\binits{H.~W.}},
  \bauthor{\bsnm{Guerrero},~\bfnm{R.~A.}\binits{R.~A.}},
  \bauthor{\bsnm{Favero},~\bfnm{M.}\binits{M.}} \AND
  \bauthor{\bsnm{Bava},~\bfnm{J.}\binits{J.}}
(\byear{2004}).
\btitle{Marine fronts at the continental shelves of austral South America:
  Physical and ecological processes}.
\bjournal{Journal of Marine Systems}
\bvolume{44}
\bpages{83--105}.%
\bptok{imsref}%
\end{barticle}%
\endbibitem

\bibitem[\protect\citeauthoryear{Appuhamillage and Sheldon}{2012}]{appsheldon}
\begin{barticle}[mr]
\bauthor{\bsnm{Appuhamillage},~\bfnm{Thilanka}\binits{T.}} \AND
  \bauthor{\bsnm{Sheldon},~\bfnm{Daniel}\binits{D.}}
(\byear{2012}).
\btitle{First passage time of skew {B}rownian motion}.
\bjournal{J. Appl. Probab.}
\bvolume{49}
\bpages{685--696}.
\bid{doi={10.1239/jap/1346955326}, issn={0021-9002}, mr={3012092}}
\bptok{imsref}%
\end{barticle}
\endbibitem

\bibitem[\protect\citeauthoryear{Appuhamillage et~al.}{2010}]{AppuhamilageWRR}
\begin{barticle}[auto:STB|2013/09/19|12:14:10]
\bauthor{\bsnm{Appuhamillage},~\bfnm{T.~A.}\binits{T.~A.}},
  \bauthor{\bsnm{Bokil},~\bfnm{V.~A.}\binits{V.~A.}},
  \bauthor{\bsnm{Thomann},~\bfnm{E.~A.}\binits{E.~A.}},
  \bauthor{\bsnm{Waymire},~\bfnm{E.~C.}\binits{E.~C.}} \AND
  \bauthor{\bsnm{Wood},~\bfnm{B.~D.}\binits{B.~D.}}
(\byear{2010}).
\btitle{Solute transport across an interface: A Fickian theory for skewness in
  breakthrough curves}.
\bjournal{Water Resour. Res.}
\bvolume{46}
\bpages{W07511}.
\bptok{imsref}%
\end{barticle}
\endbibitem

\bibitem[\protect\citeauthoryear{Appuhamillage
  et~al.}{2011a}]{AppuhamAAPCorrect}
\begin{barticle}[mr]
\bauthor{\bsnm{Appuhamillage},~\bfnm{Thilanka}\binits{T.}},
  \bauthor{\bsnm{Bokil},~\bfnm{Vrushali}\binits{V.}},
  \bauthor{\bsnm{Thomann},~\bfnm{Enrique}\binits{E.}},
  \bauthor{\bsnm{Waymire},~\bfnm{Edward}\binits{E.}} \AND
  \bauthor{\bsnm{Wood},~\bfnm{Brian}\binits{B.}}
(\byear{2011}a).
\btitle{Corrections for ``{O}ccupation and local times for skew {B}rownian
  motion with applications to dispersion across an interface.''}
\bjournal{Ann. Appl. Probab.}
\bvolume{21}
\bpages{2050--2051}.
\bid{mr={2884059}}
\bptok{imsref}%
\end{barticle}
\endbibitem

\bibitem[\protect\citeauthoryear{Appuhamillage et~al.}{2011b}]{AppuhamilageAAP}
\begin{barticle}[mr]
\bauthor{\bsnm{Appuhamillage},~\bfnm{Thilanka}\binits{T.}},
  \bauthor{\bsnm{Bokil},~\bfnm{Vrushali}\binits{V.}},
  \bauthor{\bsnm{Thomann},~\bfnm{Enrique}\binits{E.}},
  \bauthor{\bsnm{Waymire},~\bfnm{Edward}\binits{E.}} \AND
  \bauthor{\bsnm{Wood},~\bfnm{Brian}\binits{B.}}
(\byear{2011}b).
\btitle{Occupation and local times for skew {B}rownian motion with applications
  to dispersion across an interface}.
\bjournal{Ann. Appl. Probab.}
\bvolume{21}
\bpages{183--214}.
\bid{doi={10.1214/10-AAP691}, issn={1050-5164}, mr={2759199}}
\bptok{imsref}%
\end{barticle}
\endbibitem

\bibitem[\protect\citeauthoryear{Appuhamillage
  et~al.}{2012}]{Appuhamillage2012}
\begin{bmisc}[auto:STB|2013/09/19|12:14:10]
\bauthor{\bsnm{Appuhamillage},~\bfnm{T.~A.}\binits{T.~A.}},
  \bauthor{\bsnm{Bokil},~\bfnm{V.~A.}\binits{V.~A.}},
  \bauthor{\bsnm{Thomann},~\bfnm{E.~A.}\binits{E.~A.}},
  \bauthor{\bsnm{Waymire},~\bfnm{E.~C.}\binits{E.~C.}} \AND
  \bauthor{\bsnm{Wood},~\bfnm{B.~D.}\binits{B.~D.}}
(\byear{2012}).
\bhowpublished{Interfacial phenomena and natural local time. Unpublished
  manuscript. Available at \arxivurl{arXiv:1204.0271}}.
\bptok{imsref}%
\end{bmisc}
\endbibitem

\bibitem[\protect\citeauthoryear{Aris}{1956}]{aris}
\begin{barticle}[auto:STB|2013/09/19|12:14:10]
\bauthor{\bsnm{Aris},~\bfnm{R.}\binits{R.}}
(\byear{1956}).
\btitle{On the dispersion of a solute particle in a fluid moving through a
  tube}.
\bjournal{Proc. R. Soc. Lond. Ser. A Math. Phys. Eng. Sci.}
\bvolume{235}
\bpages{67--77}.
\bptok{imsref}%
\end{barticle}
\endbibitem

\bibitem[\protect\citeauthoryear{Aryasova and Portenko}{2008}]{portenko}
\begin{barticle}[mr]
\bauthor{\bsnm{Aryasova},~\bfnm{Olga~V.}\binits{O.~V.}} \AND
  \bauthor{\bsnm{Portenko},~\bfnm{Mykola~I.}\binits{M.~I.}}
(\byear{2008}).
\btitle{A uniqueness theorem for the martingale problem describing a diffusion
  in media with membranes}.
\bjournal{Theory Stoch. Process.}
\bvolume{14}
\bpages{1--9}.
\bid{issn={0321-3900}, mr={2479727}}
\bptok{imsref}%
\end{barticle}
\endbibitem

\bibitem[\protect\citeauthoryear{Barndorff-Nielsen}{1998}]{barndorff}
\begin{bincollection}[auto:STB|2013/09/19|12:14:10]
\bauthor{\bsnm{Barndorff-Nielsen},~\bfnm{O.~E.}\binits{O.~E.}}
(\byear{1998}).
\btitle{Stochastic Methods in Hydrology: Rain, Landforms and Floods, Volume~7}.
In \bbooktitle{Advanced Series on Statistical Science and Applied
Probability}.
\bpublisher{World Scientific}, \blocation{Singapore}.
\bptok{imsref}%
\end{bincollection}
\endbibitem

\bibitem[\protect\citeauthoryear{Berkowitz et~al.}{2009}]{Berkowitz09}
\begin{barticle}[auto:STB|2013/09/19|12:14:10]
\bauthor{\bsnm{Berkowitz},~\bfnm{B.}\binits{B.}},
  \bauthor{\bsnm{Cortis},~\bfnm{A.}\binits{A.}},
  \bauthor{\bsnm{Dror},~\bfnm{I.}\binits{I.}} \AND
  \bauthor{\bsnm{Scher},~\bfnm{H.}\binits{H.}}
(\byear{2009}).
\btitle{Laboratory experiments on dispersive transport across interfaces: The
  role of flow direction}.
\bjournal{Water Resour. Res.}
\bvolume{45}
\bpages{2}.
\bptok{imsref}%
\end{barticle}
\endbibitem

\bibitem[\protect\citeauthoryear{Bhattacharya and Gupta}{1984}]{bhattgupta84}
\begin{barticle}[mr]
\bauthor{\bsnm{Bhattacharya},~\bfnm{R.~N.}\binits{R.~N.}} \AND
  \bauthor{\bsnm{Gupta},~\bfnm{Vijay~K.}\binits{V.~K.}}
(\byear{1984}).
\btitle{On the {T}aylor--{A}ris theory of solute transport in a capillary}.
\bjournal{SIAM J. Appl. Math.}
\bvolume{44}
\bpages{33--39}.
\bid{doi={10.1137/0144004}, issn={0036-1399}, mr={0729999}}
\bptok{imsref}%
\end{barticle}
\endbibitem

\bibitem[\protect\citeauthoryear{Bhattacharya and Waymire}{2009}]{bhattway}
\begin{bmisc}[author]
\bauthor{\bsnm{Bhattacharya},~\bfnm{Rabi~N.}\binits{R.~N.}} \AND
  \bauthor{\bsnm{Waymire},~\bfnm{Edward~C.}\binits{E.~C.}}
(\byear{2009}).
\bhowpublished{Stochastic processes with applications.
\textit{SIAM Classics in Applied Mathematics} \textbf{61}. Originally published: Wiley, New York.}
\bid{mr={1054645}}
\bptok{imsref}%
\end{bmisc}
\endbibitem


\bibitem[\protect\citeauthoryear{Bhattacharya and Waymire}{2007}]{bcpt}
\begin{bbook}[mr]
\bauthor{\bsnm{Bhattacharya},~\bfnm{Rabi}\binits{R.}} \AND
  \bauthor{\bsnm{Waymire},~\bfnm{Edward~C.}\binits{E.~C.}}
(\byear{2007}).
\btitle{A Basic Course in Probability Theory}.
\bpublisher{Springer}, \blocation{New York}.
\bid{mr={2331066}}
\bptnote{check year}%
\bptok{imsref}%
\end{bbook}
\endbibitem

\bibitem[\protect\citeauthoryear{Bokil et~al.}{2013}]{bokil2}
\begin{bmisc}[auto:STB|2013/09/19|12:14:10]
\bauthor{\bsnm{Bokil},~\bfnm{V.~A.}\binits{V.~A.}},
  \bauthor{\bsnm{Gibson},~\bfnm{N.~L.}\binits{N.~L.}},
  \bauthor{\bsnm{Nguyen},~\bfnm{S.~T.}\binits{S.~T.}},
  \bauthor{\bsnm{Thomann},~\bfnm{E.~A.}\binits{E.~A.}} \AND
  \bauthor{\bsnm{Waymire},~\bfnm{E.~C.}\binits{E.~C.}}
(\byear{2013}).
\bhowpublished{Numerical methods for linear diffusion equations in the presence
of an interface. Preprint. Available at arXiv:\arxivurl{1310.8248}
[math.NA].}
\bptok{imsref}%
\end{bmisc}
\endbibitem

\bibitem[\protect\citeauthoryear{Brillinger}{2003}]{brill}
\begin{bincollection}[mr]
\bauthor{\bsnm{Brillinger},~\bfnm{D.~R.~R.}\binits{D.~R.~R.}}
(\byear{2003}).
\btitle{Sampling constrained animal motion using stochastic differential
  equations}.
In \bbooktitle{Probability, Statistics and Their Applications: Papers in Honor
  of Rabi Bhattacharya}.
\bseries{Institute of Mathematical Statistics Lecture Notes---Monograph Series}
\bvolume{41}
\bpages{35--48}.
\bpublisher{IMS}, \blocation{Beachwood, OH}.
\bid{mr={1999410}}
\bptok{imsref}%
\end{bincollection}
\endbibitem

\bibitem[\protect\citeauthoryear{Brillinger et~al.}{2002}]{brillinger}
\begin{barticle}[mr]
\bauthor{\bsnm{Brillinger},~\bfnm{David~R.}\binits{D.~R.}},
  \bauthor{\bsnm{Preisler},~\bfnm{Haiganoush~K.}\binits{H.~K.}},
  \bauthor{\bsnm{Ager},~\bfnm{Alan~A.}\binits{A.~A.}},
  \bauthor{\bsnm{Kie},~\bfnm{John~G.}\binits{J.~G.}} \AND
  \bauthor{\bsnm{Stewart},~\bfnm{Brent~S.}\binits{B.~S.}}
(\byear{2002}).
\btitle{Employing stochastic differential equations to model wildlife motion}.
\bjournal{Bull. Braz. Math. Soc. (N.S.)}
\bvolume{33}
\bpages{385--408}.
\bnote{Fifth Brazilian School in Probability (Ubatuba, 2001)}.
\bid{doi={10.1007/s005740200021}, issn={1678-7544}, mr={1978835}}
\bptok{imsref}%
\end{barticle}
\endbibitem

\bibitem[\protect\citeauthoryear{Brooks and Chacon}{1983}]{brookschacone}
\begin{barticle}[mr]
\bauthor{\bsnm{Brooks},~\bfnm{J.~K.}\binits{J.~K.}} \AND
  \bauthor{\bsnm{Chacon},~\bfnm{R.~V.}\binits{R.~V.}}
(\byear{1983}).
\btitle{Diffusions as a limit of stretched {B}rownian motions}.
\bjournal{Adv. Math.}
\bvolume{49}
\bpages{109--122}.
\bid{doi={10.1016/0001-8708(83)90070-1}, issn={0001-8708}, mr={0714586}}
\bptok{imsref}%
\end{barticle}
\endbibitem

\bibitem[\protect\citeauthoryear{Cantrell and
  Cosner}{2003}]{cantrell2003spatial}
\begin{bbook}[mr]
\bauthor{\bsnm{Cantrell},~\bfnm{Robert~Stephen}\binits{R.~S.}} \AND
  \bauthor{\bsnm{Cosner},~\bfnm{Chris}\binits{C.}}
(\byear{2003}).
\btitle{Spatial Ecology via Reaction--Diffusion Equations}.
\bpublisher{Wiley}, \blocation{Chichester}.
\bid{doi={10.1002/0470871296}, mr={2191264}}
\bptok{imsref}%
\end{bbook}
\endbibitem

\bibitem[\protect\citeauthoryear{Carslaw and Jaeger}{1988}]{CarslawJaeger}
\begin{bbook}[mr]
\bauthor{\bsnm{Carslaw},~\bfnm{H.~S.}\binits{H.~S.}} \AND
  \bauthor{\bsnm{Jaeger},~\bfnm{J.~C.}\binits{J.~C.}}
(\byear{1988}).
\btitle{Conduction of Heat in Solids},
\bedition{2nd} ed.
\bpublisher{Oxford Univ. Press}, \blocation{New York}.
\bid{mr={0959730}}
\bptok{imsref}%
\end{bbook}
\endbibitem

\bibitem[\protect\citeauthoryear{Chen and Fukushima}{2012}]{chenfuki}
\begin{bbook}[mr]
\bauthor{\bsnm{Chen},~\bfnm{Zhen-Qing}\binits{Z.-Q.}} \AND
  \bauthor{\bsnm{Fukushima},~\bfnm{Masatoshi}\binits{M.}}
(\byear{2012}).
\btitle{Symmetric {M}arkov Processes, Time Change, and Boundary Theory}.
\bseries{London Mathematical Society Monographs Series}
\bvolume{35}.
\bpublisher{Princeton Univ. Press}, \blocation{Princeton, NJ}.
\bid{mr={2849840}}
\bptok{imsref}%
\end{bbook}
\endbibitem

\bibitem[\protect\citeauthoryear{Cherny, Shiryaev and
  Yor}{2002}]{chernyshiryor}
\begin{barticle}[mr]
\bauthor{\bsnm{Cherny},~\bfnm{A.~S.}\binits{A.~S.}},
  \bauthor{\bsnm{Shiryaev},~\bfnm{A.~N.}\binits{A.~N.}} \AND
  \bauthor{\bsnm{Yor},~\bfnm{M.}\binits{M.}}
(\byear{2002}).
\btitle{Limit behaviour of the ``horizontal--vertical'' random walk and some
  extensions of the {D}onsker--{P}rokhorov invariance principle}.
\bjournal{Teor. Veroyatn. Primen.}
\bvolume{47}
\bpages{498--517}.
\bid{doi={10.1137/S0040585X97979834}, issn={0040-361X}, mr={1975425}}
\bptok{imsref}%
\end{barticle}
\endbibitem

\bibitem[\protect\citeauthoryear{Decamps, Goovaerts and
  Schoutens}{2006}]{decamps2006}
\begin{barticle}[mr]
\bauthor{\bsnm{Decamps},~\bfnm{Marc}\binits{M.}},
  \bauthor{\bsnm{Goovaerts},~\bfnm{Marc}\binits{M.}} \AND
  \bauthor{\bsnm{Schoutens},~\bfnm{Wim}\binits{W.}}
(\byear{2006}).
\btitle{Asymmetric skew {B}essel processes and their applications to finance}.
\bjournal{J. Comput. Appl. Math.}
\bvolume{186}
\bpages{130--147}.
\bid{doi={10.1016/j.cam.2005.03.067}, issn={0377-0427}, mr={2190302}}
\bptok{imsref}%
\end{barticle}
\endbibitem

\bibitem[\protect\citeauthoryear{{\'E}tor{\'e} and Martinez}{2013}]{Etore}
\begin{barticle}[mr]
\bauthor{\bsnm{{\'E}tor{\'e}},~\bfnm{Pierre}\binits{P.}} \AND
  \bauthor{\bsnm{Martinez},~\bfnm{Miguel}\binits{M.}}
(\byear{2013}).
\btitle{Exact simulation of one-dimensional stochastic differential equations
  involving the local time at zero of the unknown process}.
\bjournal{Monte Carlo Methods Appl.}
\bvolume{19}
\bpages{41--71}.
\bid{doi={10.1515/mcma-2013-0002}, issn={0929-9629}, mr={3039402}}
\bptok{imsref}%
\end{barticle}
\endbibitem

\bibitem[\protect\citeauthoryear{Fagan}{2002}]{fagan2}
\begin{barticle}[auto:STB|2013/09/19|12:14:10]
\bauthor{\bsnm{Fagan},~\bfnm{W.~F.}\binits{W.~F.}}
(\byear{2002}).
\btitle{Connectivity, fragmentation, and extinction risk in dendritic
  metapopulations}.
\bjournal{Ecology}
\bvolume{83}
\bpages{3243--3249}.
\bptok{imsref}%
\end{barticle}
\endbibitem

\bibitem[\protect\citeauthoryear{Fagan, Cantrell and Cosner}{1999}]{fagan1}
\begin{barticle}[auto:STB|2013/09/19|12:14:10]
\bauthor{\bsnm{Fagan},~\bfnm{W.~F.}\binits{W.~F.}},
  \bauthor{\bsnm{Cantrell},~\bfnm{R.~S.}\binits{R.~S.}} \AND
  \bauthor{\bsnm{Cosner},~\bfnm{C.}\binits{C.}}
(\byear{1999}).
\btitle{How habitat edges change species interactions}.
\bjournal{The American Naturalist}
\bvolume{153}
\bpages{165--182}.
\bptok{imsref}%
\end{barticle}
\endbibitem

\bibitem[\protect\citeauthoryear{Felder and Waymire}{2013}]{felderway}
\begin{bmisc}[auto:STB|2013/09/19|12:14:10]
\bauthor{\bsnm{Felder},~\bfnm{W.}\binits{W.}} \AND
  \bauthor{\bsnm{Waymire},~\bfnm{E.~C.}\binits{E.~C.}}
(\byear{2013}).
\bhowpublished{On the drift paradox in a regime-switching mode. Preprint.
  Available at \arxivurl{arXiv:1304.1208}}.
\bptok{imsref}%
\end{bmisc}
\endbibitem

\bibitem[\protect\citeauthoryear{Fernholz, Ichiba and
  Karatzas}{2013}]{Fernholz2012fk}
\begin{barticle}[mr]
\bauthor{\bsnm{Fernholz},~\bfnm{E.~Robert}\binits{E.~R.}},
  \bauthor{\bsnm{Ichiba},~\bfnm{Tomoyuki}\binits{T.}} \AND
  \bauthor{\bsnm{Karatzas},~\bfnm{Ioannis}\binits{I.}}
(\byear{2013}).
\btitle{Two {B}rownian particles with rank-based characteristics and
  skew-elastic collisions}.
\bjournal{Stochastic Process. Appl.}
\bvolume{123}
\bpages{2999--3026}.
\bid{doi={10.1016/j.spa.2013.03.019}, issn={0304-4149}, mr={3062434}}
\bptnote{check year}%
\bptok{imsref}%
\end{barticle}
\endbibitem

\bibitem[\protect\citeauthoryear{Fernholz et~al.}{2013}]{fernholz2011planar}
\begin{barticle}[mr]
\bauthor{\bsnm{Fernholz},~\bfnm{E.~Robert}\binits{E.~R.}},
  \bauthor{\bsnm{Ichiba},~\bfnm{Tomoyuki}\binits{T.}},
  \bauthor{\bsnm{Karatzas},~\bfnm{Ioannis}\binits{I.}} \AND
  \bauthor{\bsnm{Prokaj},~\bfnm{Vilmos}\binits{V.}}
(\byear{2013}).
\btitle{Planar diffusions with rank-based characteristics and perturbed
  {T}anaka equations}.
\bjournal{Probab. Theory Related Fields}
\bvolume{156}
\bpages{343--374}.
\bid{doi={10.1007/s00440-012-0430-7}, issn={0178-8051}, mr={3055262}}
\bptnote{check year}%
\bptok{imsref}%
\end{barticle}
\endbibitem

\bibitem[\protect\citeauthoryear{Freidlin and Sheu}{2000}]{Freidlin00}
\begin{barticle}[mr]
\bauthor{\bsnm{Freidlin},~\bfnm{Mark}\binits{M.}} \AND
  \bauthor{\bsnm{Sheu},~\bfnm{Shuenn-Jyi}\binits{S.-J.}}
(\byear{2000}).
\btitle{Diffusion processes on graphs: Stochastic differential equations, large
  deviation principle}.
\bjournal{Probab. Theory Related Fields}
\bvolume{116}
\bpages{181--220}.
\bid{doi={10.1007/PL00008726}, issn={0178-8051}, mr={1743769}}
\bptok{imsref}%
\end{barticle}
\endbibitem

\bibitem[\protect\citeauthoryear{Freidlin and Wentzell}{1993}]{Freidlin1993}
\begin{barticle}[mr]
\bauthor{\bsnm{Freidlin},~\bfnm{Mark~I.}\binits{M.~I.}} \AND
  \bauthor{\bsnm{Wentzell},~\bfnm{Alexander~D.}\binits{A.~D.}}
(\byear{1993}).
\btitle{Diffusion processes on graphs and the averaging principle}.
\bjournal{Ann. Probab.}
\bvolume{21}
\bpages{2215--2245}.
\bid{issn={0091-1798}, mr={1245308}}
\bptok{imsref}%
\end{barticle}
\endbibitem

\bibitem[\protect\citeauthoryear{Fukushima, {\=O}shima and
  Takeda}{1994}]{Fukushima94}
\begin{bbook}[mr]
\bauthor{\bsnm{Fukushima},~\bfnm{Masatoshi}\binits{M.}},
  \bauthor{\bsnm{{\=O}shima},~\bfnm{Y{\=o}ichi}\binits{Y.}} \AND
  \bauthor{\bsnm{Takeda},~\bfnm{Masayoshi}\binits{M.}}
(\byear{1994}).
\btitle{Dirichlet Forms and Symmetric {M}arkov Processes}.
\bseries{de Gruyter Studies in Mathematics}
\bvolume{19}.
\bpublisher{de Gruyter}, \blocation{Berlin}.
\bid{doi={10.1515/9783110889741}, mr={1303354}}
\bptok{imsref}%
\end{bbook}
\endbibitem

\bibitem[\protect\citeauthoryear{Gelhar and Axness}{1983}]{gelharaxness1983}
\begin{barticle}[auto:STB|2013/09/19|12:14:10]
\bauthor{\bsnm{Gelhar},~\bfnm{L.~W.}\binits{L.~W.}} \AND
  \bauthor{\bsnm{Axness},~\bfnm{C.~L.}\binits{C.~L.}}
(\byear{1983}).
\btitle{Three-dimensional stochastic analysis of macrodispersion in aquifers}.
\bjournal{Water Resour. Res.}
\bvolume{19}
\bpages{161--180}.
\bptok{imsref}%
\end{barticle}
\endbibitem

\bibitem[\protect\citeauthoryear{Gutierrez et~al.}{2012}]{Gutierrez2011be}
\begin{barticle}[mr]
\bauthor{\bsnm{Gutierrez},~\bfnm{Juan~B.}\binits{J.~B.}},
  \bauthor{\bsnm{Hurdal},~\bfnm{Monica~K.}\binits{M.~K.}},
  \bauthor{\bsnm{Parshad},~\bfnm{Rana~D.}\binits{R.~D.}} \AND
  \bauthor{\bsnm{Teem},~\bfnm{John~L.}\binits{J.~L.}}
(\byear{2012}).
\btitle{Analysis of the {T}rojan {Y} chromosome model for eradication of
  invasive species in a dendritic riverine system}.
\bjournal{J. Math. Biol.}
\bvolume{64}
\bpages{319--340}.
\bid{doi={10.1007/s00285-011-0413-9}, issn={0303-6812}, mr={2864846}}
\bptnote{check year}%
\bptok{imsref}%
\end{barticle}
\endbibitem

\bibitem[\protect\citeauthoryear{Hairer and Manson}{2010}]{hairer1}
\begin{barticle}[mr]
\bauthor{\bsnm{Hairer},~\bfnm{Martin}\binits{M.}} \AND
  \bauthor{\bsnm{Manson},~\bfnm{Charles}\binits{C.}}
(\byear{2010}).
\btitle{Periodic homogenization with an interface: The one-dimensional case}.
\bjournal{Stochastic Process. Appl.}
\bvolume{120}
\bpages{1589--1605}.
\bid{doi={10.1016/j.spa.2010.03.016}, issn={0304-4149}, mr={2653267}}
\bptok{imsref}%
\end{barticle}
\endbibitem

\bibitem[\protect\citeauthoryear{Hairer and Manson}{2011}]{hairer2}
\begin{barticle}[mr]
\bauthor{\bsnm{Hairer},~\bfnm{Martin}\binits{M.}} \AND
  \bauthor{\bsnm{Manson},~\bfnm{Charles}\binits{C.}}
(\byear{2011}).
\btitle{Periodic homogenization with an interface: The multi-dimensional case}.
\bjournal{Ann. Probab.}
\bvolume{39}
\bpages{648--682}.
\bid{doi={10.1214/10-AOP564}, issn={0091-1798}, mr={2789509}}
\bptok{imsref}%
\end{barticle}
\endbibitem

\bibitem[\protect\citeauthoryear{Harrison and Shepp}{1981}]{harrisonshepp}
\begin{barticle}[mr]
\bauthor{\bsnm{Harrison},~\bfnm{J.~M.}\binits{J.~M.}} \AND
  \bauthor{\bsnm{Shepp},~\bfnm{L.~A.}\binits{L.~A.}}
(\byear{1981}).
\btitle{On skew {B}rownian motion}.
\bjournal{Ann. Probab.}
\bvolume{9}
\bpages{309--313}.
\bid{issn={0091-1798}, mr={0606993}}
\bptok{imsref}%
\end{barticle}
\endbibitem

\bibitem[\protect\citeauthoryear{Hoteit et~al.}{2002}]{hoteit}
\begin{barticle}[mr]
\bauthor{\bsnm{Hoteit},~\bfnm{H.}\binits{H.}},
  \bauthor{\bsnm{Mose},~\bfnm{R.}\binits{R.}},
  \bauthor{\bsnm{Younes},~\bfnm{A.}\binits{A.}},
  \bauthor{\bsnm{Lehmann},~\bfnm{F.}\binits{F.}} \AND
  \bauthor{\bsnm{Ackerer},~\bfnm{Ph.}\binits{P.}}
(\byear{2002}).
\btitle{Three-dimensional modeling of mass transfer in porous media using the
  mixed hybrid finite elements and the random-walk methods}.
\bjournal{Math. Geol.}
\bvolume{34}
\bpages{435--456}.
\bid{doi={10.1023/A:1015083111971}, issn={0882-8121}, mr={1951790}}
\bptok{imsref}%
\end{barticle}
\endbibitem

\bibitem[\protect\citeauthoryear{It{\^o} and McKean}{1963}]{IM1963}
\begin{barticle}[mr]
\bauthor{\bsnm{It{\^o}},~\bfnm{K.}\binits{K.}} \AND
  \bauthor{\bsnm{McKean},~\bfnm{H.~P.}\binits{H.~P.} \bsuffix{Jr.}}
(\byear{1963}).
\btitle{Brownian motions on a half line}.
\bjournal{Illinois J. Math.}
\bvolume{7}
\bpages{181--231}.
\bid{issn={0019-2082}, mr={0154338}}
\bptok{imsref}%
\end{barticle}
\endbibitem

\bibitem[\protect\citeauthoryear{Karatzas and Shreve}{1988}]{karshrevebook}
\begin{bbook}[mr]
\bauthor{\bsnm{Karatzas},~\bfnm{Ioannis}\binits{I.}} \AND
  \bauthor{\bsnm{Shreve},~\bfnm{Steven~E.}\binits{S.~E.}}
(\byear{1988}).
\btitle{Brownian Motion and Stochastic Calculus}.
\bseries{Graduate Texts in Mathematics}
\bvolume{113}.
\bpublisher{Springer}, \blocation{New York}.
\bid{doi={10.1007/978-1-4684-0302-2}, mr={0917065}}
\bptok{imsref}%
\end{bbook}
\endbibitem

\bibitem[\protect\citeauthoryear{Kuo et~al.}{1999}]{kuo1999}
\begin{barticle}[auto:STB|2013/09/19|12:14:10]
\bauthor{\bsnm{Kuo},~\bfnm{R.~H.}\binits{R.~H.}},
  \bauthor{\bsnm{Irwin},~\bfnm{N.~C.}\binits{N.~C.}},
  \bauthor{\bsnm{Greenkorn},~\bfnm{R.~A.}\binits{R.~A.}} \AND
  \bauthor{\bsnm{Cushman},~\bfnm{J.~H.}\binits{J.~H.}}
(\byear{1999}).
\btitle{Experimental investigation of mixing in aperiodic heterogeneous porous
  media: Comparison with stochastic transport theory}.
\bjournal{Transp. Porous Media}
\bvolume{37}
\bpages{169--182}.
\bptok{imsref}%
\end{barticle}
\endbibitem

\bibitem[\protect\citeauthoryear{LaBolle, Quastel and Fogg}{1998}]{labolle}
\begin{barticle}[auto:STB|2013/09/19|12:14:10]
\bauthor{\bsnm{LaBolle},~\bfnm{E.~M.}\binits{E.~M.}},
  \bauthor{\bsnm{Quastel},~\bfnm{J.}\binits{J.}} \AND
  \bauthor{\bsnm{Fogg},~\bfnm{G.~E.}\binits{G.~E.}}
(\byear{1998}).
\btitle{Diffusion theory for transport in porous media: Transition-probability
  densities of diffusion processes corresponding to advection--dispersion
  equations}.
\bjournal{Water Resour. Res.}
\bvolume{34}
\bpages{1685--1693}.
\bptok{imsref}%
\end{barticle}
\endbibitem

\bibitem[\protect\citeauthoryear{LaBolle et~al.}{2000}]{gravner}
\begin{barticle}[auto:STB|2013/09/19|12:14:10]
\bauthor{\bsnm{LaBolle},~\bfnm{E.~M.}\binits{E.~M.}},
  \bauthor{\bsnm{Quastel},~\bfnm{J.}\binits{J.}},
  \bauthor{\bsnm{Fogg},~\bfnm{G.~E.}\binits{G.~E.}} \AND
  \bauthor{\bsnm{Gravner},~\bfnm{J.}\binits{J.}}
(\byear{2000}).
\btitle{Diffusion processes in composite porous media and their numerical
  integration by random walks: Generalized stochastic differential equations
  with discontinuous coefficients}.
\bjournal{Water Resour. Res.}
\bvolume{36}
\bpages{1685--1693}.
\bptok{imsref}%
\end{barticle}
\endbibitem

\bibitem[\protect\citeauthoryear{Le~Gall}{1984}]{LeGall}
\begin{bincollection}[mr]
\bauthor{\bsnm{Le~Gall},~\bfnm{J.~F.}\binits{J.~F.}}
(\byear{1984}).
\btitle{One-dimensional stochastic differential equations involving the local
  times of the unknown process}.
In \bbooktitle{Stochastic Analysis and Applications ({S}wansea, 1983)}.
\bseries{Lecture Notes in Math.}
\bvolume{1095}
\bpages{51--82}.
\bpublisher{Springer}, \blocation{Berlin}.
\bid{doi={10.1007/BFb0099122}, mr={0777514}}
\bptok{imsref}%
\end{bincollection}
\endbibitem

\bibitem[\protect\citeauthoryear{Lejay}{2006}]{lejaysurvey}
\begin{barticle}[mr]
\bauthor{\bsnm{Lejay},~\bfnm{Antoine}\binits{A.}}
(\byear{2006}).
\btitle{On the constructions of the skew {B}rownian motion}.
\bjournal{Probab. Surv.}
\bvolume{3}
\bpages{413--466}.
\bid{doi={10.1214/154957807000000013}, issn={1549-5787}, mr={2280299}}
\bptok{imsref}%
\end{barticle}
\endbibitem

\bibitem[\protect\citeauthoryear{Lejay and Pichot}{2012}]{lejay2012simulating}
\begin{barticle}[mr]
\bauthor{\bsnm{Lejay},~\bfnm{Antoine}\binits{A.}} \AND
  \bauthor{\bsnm{Pichot},~\bfnm{G{\'e}raldine}\binits{G.}}
(\byear{2012}).
\btitle{Simulating diffusion processes in discontinuous media: A numerical
  scheme with constant time steps}.
\bjournal{J. Comput. Phys.}
\bvolume{231}
\bpages{7299--7314}.
\bid{doi={10.1016/j.jcp.2012.07.011}, issn={0021-9991}, mr={2969713}}
\bptok{imsref}%
\end{barticle}
\endbibitem

\bibitem[\protect\citeauthoryear{Lutscher, Lewis and
  McCauley}{2006}]{lutscheretal2006}
\begin{barticle}[mr]
\bauthor{\bsnm{Lutscher},~\bfnm{Frithjof}\binits{F.}},
  \bauthor{\bsnm{Lewis},~\bfnm{Mark~A.}\binits{M.~A.}} \AND
  \bauthor{\bsnm{McCauley},~\bfnm{Edward}\binits{E.}}
(\byear{2006}).
\btitle{Effects of heterogeneity on spread and persistence in rivers}.
\bjournal{Bull. Math. Biol.}
\bvolume{68}
\bpages{2129--2160}.
\bid{doi={10.1007/s11538-006-9100-1}, issn={0092-8240}, mr={2293837}}
\bptok{imsref}%
\end{barticle}
\endbibitem

\bibitem[\protect\citeauthoryear{Lutscher, Pachepsky and
  Lewis}{2005}]{lutscheretal2005}
\begin{barticle}[mr]
\bauthor{\bsnm{Lutscher},~\bfnm{Frithjof}\binits{F.}},
  \bauthor{\bsnm{Pachepsky},~\bfnm{Elizaveta}\binits{E.}} \AND
  \bauthor{\bsnm{Lewis},~\bfnm{Mark~A.}\binits{M.~A.}}
(\byear{2005}).
\btitle{The effect of dispersal patterns on stream populations}.
\bjournal{SIAM Rev.}
\bvolume{47}
\bpages{749--772 (electronic)}.
\bid{doi={10.1137/050636152}, issn={0036-1445}, mr={2212398}}
\bptok{imsref}%
\end{barticle}
\endbibitem

\bibitem[\protect\citeauthoryear{Ma and R{\"o}ckner}{1992}]{marockner}
\begin{bbook}[mr]
\bauthor{\bsnm{Ma},~\bfnm{Zhi~Ming}\binits{Z.~M.}} \AND
  \bauthor{\bsnm{R{\"o}ckner},~\bfnm{Michael}\binits{M.}}
(\byear{1992}).
\btitle{Introduction to the Theory of (Nonsymmetric) {D}irichlet Forms}.
\bpublisher{Springer}, \blocation{Berlin}.
\bid{doi={10.1007/978-3-642-77739-4}, mr={1214375}}
\bptok{imsref}%
\end{bbook}
\endbibitem

\bibitem[\protect\citeauthoryear{Martinez and Talay}{2012}]{martineztalay}
\begin{barticle}[mr]
\bauthor{\bsnm{Martinez},~\bfnm{Miguel}\binits{M.}} \AND
  \bauthor{\bsnm{Talay},~\bfnm{Denis}\binits{D.}}
(\byear{2012}).
\btitle{One-dimensional parabolic diffraction equations: Pointwise estimates
  and discretization of related stochastic differential equations with weighted
  local times}.
\bjournal{Electron. J. Probab.}
\bvolume{17}
\bpages{1--32}.
\bid{doi={10.1214/EJP.v17-1905}, issn={1083-6489}, mr={2912504}}
\bptok{imsref}%
\end{barticle}
\endbibitem

\bibitem[\protect\citeauthoryear{Matano and Palma}{2008}]{MatanoPalma}
\begin{barticle}[auto:STB|2013/09/19|12:14:10]
\bauthor{\bsnm{Matano},~\bfnm{R.}\binits{R.}} \AND
  \bauthor{\bsnm{Palma},~\bfnm{E.}\binits{E.}}
(\byear{2008}).
\btitle{On the upwelling of downwelling currents}.
\bjournal{Journal of Physical Oceanography}
\bvolume{38}
\bpages{2482--2500}.
\bptok{imsref}%
\end{barticle}
\endbibitem

\bibitem[\protect\citeauthoryear{M{\"u}ller}{1954}]{muller1954investigations}
\begin{barticle}[auto:STB|2013/09/19|12:14:10]
\bauthor{\bsnm{M{\"u}ller},~\bfnm{K.}\binits{K.}}
(\byear{1954}).
\btitle{Investigations on the organic drift in north swedish streams}.
\bjournal{Rep. Inst. Freshwat. Res. Drottningholm}
\bvolume{35}
\bpages{133--148}.
\bptok{imsref}%
\end{barticle}
\endbibitem

\bibitem[\protect\citeauthoryear{Nakao}{1972}]{Nakao}
\begin{barticle}[mr]
\bauthor{\bsnm{Nakao},~\bfnm{Shintaro}\binits{S.}}
(\byear{1972}).
\btitle{On the pathwise uniqueness of solutions of one-dimensional stochastic
  differential equations}.
\bjournal{Osaka J. Math.}
\bvolume{9}
\bpages{513--518}.
\bid{issn={0030-6126}, mr={0326840}}
\bptok{imsref}%
\end{barticle}
\endbibitem

\bibitem[\protect\citeauthoryear{Okubo and Levin}{2001}]{levin}
\begin{bbook}[mr]
\bauthor{\bsnm{Okubo},~\bfnm{Akira}\binits{A.}} \AND
  \bauthor{\bsnm{Levin},~\bfnm{Simon~A.}\binits{S.~A.}}
(\byear{2001}).
\btitle{Diffusion and Ecological Problems: Modern Perspectives},
\bedition{2nd} ed.
\bseries{Interdisciplinary Applied Mathematics}
\bvolume{14}.
\bpublisher{Springer}, \blocation{New York}.
\bid{mr={1895041}}
\bptnote{check year}%
\bptok{imsref}%
\end{bbook}
\endbibitem

\bibitem[\protect\citeauthoryear{Ouknine}{1990}]{Ouknine}
\begin{barticle}[mr]
\bauthor{\bsnm{Ouknine},~\bfnm{Y.}\binits{Y.}}
(\byear{1990}).
\btitle{Le ``{S}kew-{B}rownian motion'' et les processus qui en d\'erivent}.
\bjournal{Teor. Veroyatn. Primen.}
\bvolume{35}
\bpages{173--179}.
\bid{doi={10.1137/1135018}, issn={0040-361X}, mr={1050069}}
\bptok{imsref}%
\end{barticle}
\endbibitem

\bibitem[\protect\citeauthoryear{Ovaskainen and Cornell}{2003}]{ovas}
\begin{barticle}[mr]
\bauthor{\bsnm{Ovaskainen},~\bfnm{Otso}\binits{O.}} \AND
  \bauthor{\bsnm{Cornell},~\bfnm{Stephen~J.}\binits{S.~J.}}
(\byear{2003}).
\btitle{Biased movement at a boundary and conditional occupancy times for
  diffusion processes}.
\bjournal{J. Appl. Probab.}
\bvolume{40}
\bpages{557--580}.
\bid{issn={0021-9002}, mr={1993253}}
\bptok{imsref}%
\end{barticle}
\endbibitem

\bibitem[\protect\citeauthoryear{Peckham}{1995}]{peckham}
\begin{barticle}[auto:STB|2013/09/19|12:14:10]
\bauthor{\bsnm{Peckham},~\bfnm{S.~D.}\binits{S.~D.}}
(\byear{1995}).
\btitle{New results for self-similar trees with applications to river
  networks}.
\bjournal{Water Resources Research}
\bvolume{31}
\bpages{1023--1029}.
\bptok{imsref}%
\end{barticle}
\endbibitem

\bibitem[\protect\citeauthoryear{Perrin}{1913}]{perrin}
\begin{bbook}[auto:STB|2013/09/19|12:14:10]
\bauthor{\bsnm{Perrin},~\bfnm{J.}\binits{J.}}
(\byear{1913}).
\btitle{Les Atomes/Par Jean Perrin, Avec 13 Figures}.
\bpublisher{F. Alcan}, \blocation{Paris}.
\bptok{imsref}%
\end{bbook}
\endbibitem

\bibitem[\protect\citeauthoryear{Peskir}{2007}]{Peskir2007}
\begin{bincollection}[mr]
\bauthor{\bsnm{Peskir},~\bfnm{Goran}\binits{G.}}
(\byear{2007}).
\btitle{A change-of-variable formula with local time on surfaces}.
In \bbooktitle{S\'eminaire de {P}robabilit\'es {XL}}.
\bseries{Lecture Notes in Math.}
\bvolume{1899}
\bpages{69--96}.
\bpublisher{Springer}, \blocation{Berlin}.
\bid{mr={2408999}}
\bptok{imsref}%
\end{bincollection}
\endbibitem

\bibitem[\protect\citeauthoryear{Portenko}{1990}]{portenko1}
\begin{bbook}[mr]
\bauthor{\bsnm{Portenko},~\bfnm{N.~I.}\binits{N.~I.}}
(\byear{1990}).
\btitle{Generalized Diffusion Processes}.
\bseries{Translations of Mathematical Monographs}
\bvolume{83}.
\bpublisher{Amer. Math. Soc.}, \blocation{Providence, RI}.
\bid{mr={1104660}}
\bptok{imsref}%
\end{bbook}
\endbibitem

\bibitem[\protect\citeauthoryear{Prokaj}{2011}]{Prokaj}
\begin{bmisc}[auto:STB|2013/09/19|12:14:10]
\bauthor{\bsnm{Prokaj},~\bfnm{V.}\binits{V.}}
(\byear{2011}).
\bhowpublished{The solution of the perturbed Tanaka---Equation is pathwise
  unique. Available at \arxivurl{arXiv:1104.0740}}.
\bptok{imsref}%
\end{bmisc}
\endbibitem

\bibitem[\protect\citeauthoryear{Ramirez}{2011}]{ramirez2011multi}
\begin{barticle}[mr]
\bauthor{\bsnm{Ramirez},~\bfnm{Jorge~M.}\binits{J.~M.}}
(\byear{2011}).
\btitle{Multi-skewed {B}rownian motion and diffusion in layered media}.
\bjournal{Proc. Amer. Math. Soc.}
\bvolume{139}
\bpages{3739--3752}.
\bid{doi={10.1090/S0002-9939-2011-10766-4}, issn={0002-9939}, mr={2813404}}
\bptok{imsref}%
\end{barticle}
\endbibitem

\bibitem[\protect\citeauthoryear{Ramirez}{2012a}]{ramirez2011population}
\begin{barticle}[mr]
\bauthor{\bsnm{Ramirez},~\bfnm{Jorge~M.}\binits{J.~M.}}
(\byear{2012}a).
\btitle{Population persistence under advection--diffusion in river networks}.
\bjournal{J. Math. Biol.}
\bvolume{65}
\bpages{919--942}.
\bid{doi={10.1007/s00285-011-0485-6}, issn={0303-6812}, mr={2984129}}
\bptnote{check year}%
\bptok{imsref}%
\end{barticle}
\endbibitem

\bibitem[\protect\citeauthoryear{Ramirez}{2012b}]{ramirez2012green}
\begin{barticle}[mr]
\bauthor{\bsnm{Ramirez},~\bfnm{Jorge~M.}\binits{J.~M.}}
(\byear{2012}b).
\btitle{Green's functions for {S}turm--{L}iouville problems on directed tree
  graphs}.
\bjournal{Rev. Colombiana Mat.}
\bvolume{46}
\bpages{15--25}.
\bid{issn={0034-7426}, mr={2945668}}
\bptok{imsref}%
\end{barticle}
\endbibitem

\bibitem[\protect\citeauthoryear{Ramirez et~al.}{2006}]{ramirez2006}
\begin{barticle}[mr]
\bauthor{\bsnm{Ramirez},~\bfnm{Jorge~M.}\binits{J.~M.}},
  \bauthor{\bsnm{Thomann},~\bfnm{Enrique~A.}\binits{E.~A.}},
  \bauthor{\bsnm{Waymire},~\bfnm{Edward~C.}\binits{E.~C.}},
  \bauthor{\bsnm{Haggerty},~\bfnm{Roy}\binits{R.}} \AND
  \bauthor{\bsnm{Wood},~\bfnm{Brian}\binits{B.}}
(\byear{2006}).
\btitle{A generalized {T}aylor--{A}ris formula and skew diffusion}.
\bjournal{Multiscale Model. Simul.}
\bvolume{5}
\bpages{786--801}.
\bid{doi={10.1137/050642770}, issn={1540-3459}, mr={2257235}}
\bptok{imsref}%
\end{barticle}
\endbibitem

\bibitem[\protect\citeauthoryear{Revuz and Yor}{1999}]{revuzyor}
\begin{bbook}[mr]
\bauthor{\bsnm{Revuz},~\bfnm{Daniel}\binits{D.}} \AND
  \bauthor{\bsnm{Yor},~\bfnm{Marc}\binits{M.}}
(\byear{1999}).
\btitle{Continuous Martingales and {B}rownian Motion},
\bedition{3rd} ed.
\bseries{Grundlehren der Mathematischen Wissenschaften [Fundamental Principles
  of Mathematical Sciences]}
\bvolume{293}.
\bpublisher{Springer}, \blocation{Berlin}.
\bid{mr={1725357}}
\bptok{imsref}%
\end{bbook}
\endbibitem

\bibitem[\protect\citeauthoryear{Rodriguez-Iturbe and
  Rinaldo}{2001}]{RodriguezIturbe2001}
\begin{bbook}[auto:STB|2013/09/19|12:14:10]
\bauthor{\bsnm{Rodriguez-Iturbe},~\bfnm{I.}\binits{I.}} \AND
  \bauthor{\bsnm{Rinaldo},~\bfnm{A.}\binits{A.}}
(\byear{2001}).
\btitle{Fractal River Basins: Chance and Self-Organization}.
\bpublisher{Cambridge Univ. Press}, \blocation{Cambridge}.
\bptok{imsref}%
\end{bbook}
\endbibitem

\bibitem[\protect\citeauthoryear{Schultz and Crone}{2001}]{Schultz2001}
\begin{barticle}[auto:STB|2013/09/19|12:14:10]
\bauthor{\bsnm{Schultz},~\bfnm{C.~B.}\binits{C.~B.}} \AND
  \bauthor{\bsnm{Crone},~\bfnm{E.~E.}\binits{E.~E.}}
(\byear{2001}).
\btitle{Edge-mediated dispersal behavior in a prairie butterfly}.
\bjournal{Ecology}
\bvolume{82}
\bpages{1879--1892}.
\bptok{imsref}%
\end{barticle}
\endbibitem

\bibitem[\protect\citeauthoryear{Speirs and Gurney}{2001}]{Speirs2001kx}
\begin{barticle}[auto:STB|2013/09/19|12:14:10]
\bauthor{\bsnm{Speirs},~\bfnm{D.~C.}\binits{D.~C.}} \AND
  \bauthor{\bsnm{Gurney},~\bfnm{W.~S.~C.}\binits{W.~S.~C.}}
(\byear{2001}).
\btitle{Population persistence in rivers and estuaries}.
\bjournal{Ecology}
\bvolume{82}
\bpages{1219--1237}.
\bptok{imsref}%
\end{barticle}
\endbibitem

\bibitem[\protect\citeauthoryear{Stroock and Varadhan}{2006}]{StrookVaradhan}
\begin{bbook}[mr]
\bauthor{\bsnm{Stroock},~\bfnm{Daniel~W.}\binits{D.~W.}} \AND
  \bauthor{\bsnm{Varadhan},~\bfnm{S.~R.~Srinivasa}\binits{S.~R.~S.}}
(\byear{2006}).
\btitle{Multidimensional Diffusion Processes}.
\bpublisher{Springer}, \blocation{Berlin}.
\bnote{Reprint of the 1997 edition.}
\bid{mr={2190038}}
\bptnote{check year}%
\bptok{imsref}%
\end{bbook}
\endbibitem

\bibitem[\protect\citeauthoryear{Taylor}{1953}]{taylor}
\begin{barticle}[auto:STB|2013/09/19|12:14:10]
\bauthor{\bsnm{Taylor},~\bfnm{G.~I.}\binits{G.~I.}}
(\byear{1953}).
\btitle{Dispersion of a soluble matter in solvent flowing through a tube}.
\bjournal{Proc. R. Soc. Lond. Ser. A Math. Phys. Eng. Sci.}
\bvolume{219}
\bpages{186--203}.
\bptok{imsref}%
\end{barticle}
\endbibitem

\bibitem[\protect\citeauthoryear{Walsh}{1978}]{walsh1}
\begin{barticle}[auto:STB|2013/09/19|12:14:10]
\bauthor{\bsnm{Walsh},~\bfnm{J.~B.}\binits{J.~B.}}
(\byear{1978}).
\btitle{A diffusion with a discontinuous local time}.
\bjournal{Ast\'erisque}
\bvolume{52--53}
\bpages{37--45}.
\bptok{imsref}%
\end{barticle}
\endbibitem

\bibitem[\protect\citeauthoryear{Wooding}{1960}]{wooding}
\begin{barticle}[mr]
\bauthor{\bsnm{Wooding},~\bfnm{R.~A.}\binits{R.~A.}}
(\byear{1960}).
\btitle{Rayleigh instability of a thermal boundary layer in flow through a
  porous medium}.
\bjournal{J. Fluid Mech.}
\bvolume{9}
\bpages{183--192}.
\bid{issn={0022-1120}, mr={0120977}}
\bptok{imsref}%
\end{barticle}
\endbibitem

\end{thebibliography}
\end{document}